\documentclass[12pt]{amsart}
\usepackage{amsmath}
\usepackage{amssymb}
\usepackage{amscd}
\usepackage{graphicx}
\input xy
\xyoption{all}

\newtheorem{thm}{Theorem}[section]
\newtheorem{lem}[thm]{Lemma}

\newtheorem{prop}[thm]{Proposition}

\newtheorem{prob}[thm]{Problem}
\newtheorem{cor}[thm]{Corollary}

\theoremstyle{definition}
\newtheorem{rem}[thm]{Remark}
\newtheorem{de}[thm]{Definition}

\def \ol {\overline}
\def \N {\mathbb N}

\def \Z {\mathbb Z}

\def \A {\mathcal A}

\def \F {\mathcal F}

\def \P {\mathcal P}

\def \b {\beta}
\def \ep {\epsilon}

\def \lra {\longrightarrow}

\newcommand{\rest}{\upharpoonright}

\begin{document}

\title{Sufficient conditions under which a transitive system is chaotic}
\author{E. Akin, E. Glasner, W. Huang,  S. Shao and X. Ye}

\address{Mathematics Department ,
The City College, 137 Street and Convent Avenue,
New York City, NY 10031, USA}
\email{ethanakin@earthlink.net}

\address{Department of Mathematics, Tel Aviv University, Tel Aviv,
Israel}

\email{glasner@math.tau.ac.il}

\address{Department of Mathematics,
University of Science and Technology of China,
Hefei, Anhui, 230026, P.R. China}

\email{wenh@mail.ustc.edu.cn}
\email{songshao@ustc.edu.cn}
\email{yexd@ustc.edu.cn}







\begin{abstract}
Let $(X,T)$ be a topologically transitive dynamical system. We show that if
there is a subsystem $(Y,T)$ of $(X,T)$ such that
$(X\times Y, T\times T)$ is transitive, then $(X,T)$ is strongly
chaotic in the sense of Li and Yorke.
We then show that many of the known
sufficient conditions in the literature, as well as a few new results,
are corollaries of this fact.
In fact the kind of chaotic behavior we deduce
in these results is a much stronger variant of Li-Yorke chaos which we call
uniform chaos.
For minimal systems we show, among other results, that
uniform chaos is preserved by extensions and that a minimal
system which is not uniformly chaotic is PI.
 \end{abstract}

\maketitle

\tableofcontents

\section*{Introduction}

The presence or the lack of chaotic behavior is one of the most
prominent traits of a dynamical system. However, by now there exists
in the literature on dynamical systems a plethora of ways to define Chaos.
In 1975, Li and Yorke introduced a notion of chaos \cite{LY}, known now as
Li-Yorke chaos,  for interval maps. With a small modification this notion can be
extended to any metric space. Another notion was introduced later by
Devaney \cite{D}.
In \cite{GW93} the authors suggested to base the definition of chaotic behavior
on the notion of positive topological entropy.
More recently it was shown that both Devaney chaos \cite{HY02}, and
positive entropy \cite{BGKM} imply Li-Yorke chaos.
We remark that weak mixing as well (or even scattering) implies
Li-Yorke chaos. Thus, in a certain sense Li-Yorke chaos is the weakest
notion of Chaos.  We refer the reader to the recent
monograph \cite{AAG} and the review \cite{GY}
on local entropy theory, which include discussions of the
above notions.

\medskip

It is natural to ask which transitive systems are chaotic
and this is the main theme of this work.
In Section 1 we introduce our terminology and review some basic facts.
In Section 2 we first prove, the somewhat surprising fact (Theorem
\ref{tran}) that every transitive system is partially rigid. This is then used
in Section 3 to deduce the following criterion.
For a transitive topological dynamical system $(X,T)$
if there is a subsystem $(Y,T)$ of $(X,T)$
(i.e. $Y $ is a non-empty closed and $T$-invariant
subset of $X$) such that $(X\times Y, T\times T)$ is transitive, then $(X,T)$ is
strongly Li-Yorke chaotic. As we will see many of the known
sufficient conditions in the literature, as well as a few new results,
are corollaries of this fact. In fact the kind of chaotic behavior we deduce
in these results is a much stronger variant of Li-Yorke chaos which we call
uniform chaos.
In Section 4 we reexamine these results in view of the Kuratowski-Mycielski
theory.
In Section \ref{Sec-min} we specialize to minimal dynamical systems.
After reviewing some structure theory we show, among other results that
for minimal systems uniform chaos is preserved by extensions, and that
if a minimal system is not uniformly chaotic then it is a PI system.
We also show that a minimal {\em strictly PI} system which is not point distal
admits a proximal scrambled Mycielski set.
This perhaps suggests that a minimal system which does not
contain such a set is actually point distal, but we have to leave that issue
as an open problem.

Throughout the paper and mostly in Section \ref{Sec-min} we make heavy use
of enveloping semigroups and structure theory. We refer, for example,
to the sources \cite{G76}, \cite{V77},  \cite{Au}, and \cite{Ak97} for the 
necessary background.

Eli Glasner thanks his coauthors and gracious hosts, Xiangdong Ye and
Song Shao for their hospitality during a long visit to Hefei in 2004, where most
of this work was done.

\section{Preliminary definitions and results}

In this section we briefly review some basic definitions and results
from topological dynamics. Relevant references are  \cite{GM},
\cite{GW93}, \cite{AAB}, \cite{Ak97}, \cite{AG}, \cite{HY02},
\cite{Ak03}, \cite{AAG}.
The latter is perhaps a good starting point for a beginner.
One can also try to trace the historical development of these
notions from that source and the reference list thereof.

\medskip

\subsection{Transitivity and related notions} 
We write $\Z$ to denote the integers, $\Z_+$ for the non-negative
integers and $\N$ for the natural numbers. Throughout this paper a
{\em topological dynamical system} (TDS for short) is a pair
$(X,T)$, where $X$ is a non-vacuous compact metric space with a metric
$d$ and $T$ is a continuous surjective map from $X$ to itself.
A closed invariant subset $Y \subset X$ defines naturally a {\em subsystem}
$(Y,T)$ of $(X,T)$.

For subsets $A, B \subset X$ we define for a TDS $(X,T)$ the 
\emph{hitting time set}
$N(A,B):=\{n\in\Z_+: A\cap T^{-n}B\not=\emptyset\}$.
When $A=\{x\}$ is a singleton we write simply $N(x,B)$
and if moreover $B$ is a neighborhood of $x$ we refer to 
$N(x,B)$ as the set of {\em return times}.

Recall that $(X,T)$ is called {\em topologically transitive}
(or just {\em transitive}) if for every pair of nonempty open
subsets $U$ and $V$, the set $N(U,V)$ is non-empty.

Let $\omega (x,T)$ be the set of the limit points of the
orbit of $x$, 
$$
Orb(x,T):=\{x,T(x),T^2(x),\ldots\}.
$$ 
A point $x\in X$ is called a {\em transitive point} if $\omega (x,T)=X$.
It is easy to see that if $(X,T)$ is transitive then the set of
all transitive points is a dense $G_\delta$ set of $X$ (denoted by
$X_{tr}$ or $Trans(X)$). If $X_{tr}=X$ then we say that $(X,T)$ is {\em
minimal}. Equivalently, $(X,T)$ is minimal if and only if  it contains no
proper subsystems. It is well known that there is some minimal
subsystem in any dynamical system $(X,T)$, which is called a {\em
minimal set} of $X$. Each point belonging to some minimal set of
$X$ is called a {\em minimal point}.

A TDS $(X,T)$ is {\em (topologically) weakly mixing} if 
the product system$(X \times X, T \times T)$ is
transitive. 

A pair $(x,y)\in X \times X$ is said to be {\em proximal} if
$\liminf_{n \to +\infty} d(T^nx,T^ny)$ $=0$ and it is called
{\em asymptotic} when $\lim_{n \to +\infty} d(T^nx,T^ny)=0$. 
If in addition $x\neq y$, then $(x,y)$ is a {\em proper} proximal 
(or asymptotic) pair. 
The sets of proximal pairs and asymptotic
pairs of $(X,T)$ are denoted by $P(X,T)$ and $Asym(X,T)$
respectively. A point $x \in X$ is a {\em recurrent point} if there are
$n_i \nearrow +\infty$ such that $T^{n_i}x \to  x$.
A pair $(x,y) \in X^2$ which is not proximal is said to
be {\em distal}. A pair is said to be a {\em Li-Yorke pair} if it
is proximal but not asymptotic. A pair $(x,y) \in X^2 \setminus
\Delta_X $ is said to be a {\em strong Li-Yorke pair} if it is
proximal and is also a recurrent point of $X^2$. Clearly
a strong Li-Yorke pair is a Li-Yorke pair. A system without
proper proximal pairs (Li-Yorke pairs, strong Li-Yorke pairs) is
called {\em distal} ({\em almost distal}, {\em semi-distal}
respectively). It follows directly from the definitions that a distal system 
is almost distal and an almost distal system is semi-distal. 
A point $x$ is called a {\em distal point} if its proximal cell
$P[x]=\{x' \in X: (x,x') \in P(X,T)\}=\{x\}$. A system $(X,T)$
is point distal if it contains a distal point. A theorem of
Ellis \cite{E} says that in a metric minimal point distal system
the set of distal points is dense and $G_\delta$. A dynamical system
$(X,T)$ is {\em equicontinuous} if for every $\epsilon >0$ there is 
$\delta >0$ such that $d(x,y)<\delta$ implies $d(T^nx,T^ny)<\epsilon$, 
for every $n \in \Z_+$. 
Clearly an equicontinuous system is distal.

A {\em homomorphism} (or a {\em factor map}) 
$\pi : (X,T) \longrightarrow (Y,S)$ is a
continuous onto map from $X$ to $Y$ such that $S \circ \pi=\pi
\circ T$. In this situation $(X,T)$ is said to be an {\em
extension} of $(Y,S)$ and $(Y,S)$ is called a {\em factor} of $(X,T)$.
A homomorphism $\pi$ is determined by the corresponding
closed invariant equivalence relation
$R_{\pi} = \{ (x_1,x_2): \pi x_1= \pi x_2 \} =(\pi \times
\pi )^{-1} \Delta_Y \subset  X \times X$.

An extension $\pi :
(X,T) \rightarrow (Y,S)$ is called {\em asymptotic} if $R_{\pi}
\subset Asmp(X,T)$.
Similarly we define {\em proximal, distal} extensions. We define $\pi$ to be an {\em equicontiuous}
extension if for every $\epsilon >0$ there is
$\delta >0$ such that $(x,y) \in R_{\pi}$ and $d(x,y)<\delta$ implies $d(T^nx,T^ny)<\epsilon$,
for every $n \in \Z_+$. The extension $\pi$ is
called {\em almost one-to-one} if the set
$X_0=\{x \in X: \pi^{-1}(\pi(x)) = \{x\}\}$ is a
dense $G_\delta$ subset of $X$.

\medskip

\subsection{The enveloping semigroup}\label{es}


An {\em Ellis semigroup} is a semigroup equipped with a compact Hausdorff topology such that for every $p\in E$ the map 
$R_p: E \to E$ defined by $R_p(q) = qp$ is continuous. 
(This is sometimes called a right topological,
or a left topological, or a right semi-topological semigroup.
Here we try to use a non-ambiguous term which we hope will 
standardize the terminology.)
An {\em Ellis action} is an action of an Ellis semigroup $E$ on a compact Hausdorff space $X$ such that for every $x \in X$ the map $R_x : E \to X$ defined by $R_x(q) = qx$ is continuous. 

The {\em enveloping semigroup\/} $E=E(X,T)=E(X)$ of a
dynamical system $(X,T)$ is defined as the closure in $X^X$ (with
its compact, usually non-metrizable, pointwise convergence
topology) of the set $\{T^n: n \in \Z_+\}$. With the operation of
composition of maps this is an Ellis semigroup and the operation of evaluation is
an Ellis action of $E(X,T)$ on $X$ which extends the action of $\Z_+$ via $T$.


The elements of $E(X,T)$ may behave very
badly as maps of $X$ into itself; usually they are
not even Borel measurable. However our main
interest in the enveloping semigroup lies in its algebraic structure
and its dynamical significance. A key lemma
in the study of this algebraic structure is the
following:


\begin{lem}[Ellis]\label{idemp} If $E$ is an Ellis
semigroup, then $E$ contains an {\em idempotent\/};
i.e.,\ an element $v$ with $v^2=v$.
\end{lem}

In the next proposition we state some
basic properties of the enveloping semigroup
$E=E(X,T)$.

\begin{prop}\label{envel}
\begin{enumerate}
\item
A subset $I$ of $E$ is a minimal left
ideal of the semigroup $E$ if and only if it is a minimal
subsystem of $(E,T)$. In particular a
minimal left ideal is closed. We will refer to it simply
as a {\em minimal ideal\/}.
Minimal ideals $I$ in $E$ exist and for each
such ideal the set of idempotents in $I$, denoted
by $J=J(I)$, is non-empty.
\item
Let $I$ be a minimal ideal and $J$ its set
of idempotents then:
\begin{enumerate}
\item[(a)]  For $v\in J$ and $p\in I$, $pv=p$.
\item[(b)] For each $v\in J$,\ $vI=
\{vp:p\in I\}=\{p\in I: vp=p\}$ is a subgroup of $I$ with
identity element $v$.
For every $w\in J$ the map $p\mapsto wp$
is a group isomorphism of $vI$ onto $wI$.
\item[(c)] $ \{vI:v\in J\}$ is a partition
of $I$. Thus if $p\in I$ then there exists a unique
$v\in J$ such that $p\in vI$.
\end{enumerate}
\item
Let $K,L,$ and $I$ be minimal ideals of $E$.
Let $v$ be an idempotent in $I$, then there exists
a unique idempotent $v'$ in $L$ such that
$vv'=v'$ and $v'v=v$. (We write $v\sim v'$
and say that $v'$ is {\em equivalent\/} to $v$.)
If $v''\in K$ is equivalent to $v'$, then
$v''\sim v$. The map $p\mapsto pv'$ of
$I$ to $L$ is an isomorphism of dynamical systems.
\item
A pair $(x,x')\in X \times X$ is proximal
if and only if  $px=px'$ for some $p\in E$, if and only if  there exists
a minimal ideal $I$ in $E$ with $px=px'$ for
every $p\in I$.
\item
If $(X,T)$ is minimal, then the proximal cell of $x$
$$
P[x]=\{x'\in X:(x,x')\in P\}=\{vx:v\in \hat{J}\},
$$
where $\hat{J} = \bigcup \{J(I): I\ \text{is
a minimal left ideal in}\ E(X,T)\}$ is the set of
minimal idempotents.
\end{enumerate}
\end{prop}

We will make use also of the {\em adherence semigroup} $A(X,T)$
which is defined as the $\omega$-limit set of the
collection $\{T^n: n \in \Z_+\}$ in $E(X,T)$.
\medskip

Often one has to deal with more than one system at
a time; e.g., we can be working simultaneously with
a system and its factors,
two different systems, their product, subsystems of the product, etc.
Or, given a topological system $(X,T)$ we may
have to work with associated systems like the
action induced on the space $C(X)$ of closed subsets of $X$, with
its Hausdorff topology. It is therefore convenient to
have one enveloping semigroup acting on all of
the systems simultaneously. This can be easily done
by considering the enveloping semigroup of the
product of all the systems under consideration.
However, one looses nothing and gains much
in convenience as well as in added machinery if
one works instead with the ``universal" enveloping
semigroup.


Such a universal object for $\Z$-actions is $\beta\Z$, the
{\em \v{C}ech-Stone compactification} of the integers
(and $\beta\Z_+$ for $\Z_+$-actions). These are Ellis semigroups and
any $\Z$ (or $\Z_+$) action on $X$ via $T$ extends naturally to an
Ellis action of $\b \Z$ (resp. $\b \Z_+$) on $X$.



We will freely use this fact and thus will let $\beta\Z_+$  ``act" on every
compact $\Z_+$ dynamical system.
In this case the {\em corona} $\beta^*\Z_+ = \beta\Z_+ \setminus \Z_+$
coincides with the adherence semigroup.
We refer to \cite{G76}, \cite{Au}, \cite{Ak97} and \cite{G03} for more details.

\medskip

\subsection{Some notions of Chaos}

A subset $A \subset X$ is called {\em scrambled}
({\em strongly scrambled}) if every pair
of distinct points in $A$ is Li-Yorke (strong Li-Yorke).
The system $(X,T)$
is said to be {\em Li-Yorke chaotic} ({\em strong Li-Yorke} chaotic)
if it contains an uncountable scrambled (strongly scrambled) set.

The notion of equicontinuity can be localized in an obvious way.
Namely, $x \in X$ is called an {\em equicontinuity point}
if for every $\epsilon >0$ there is $\delta >0$ such that $d(x,y)<\delta$
implies $d(T^nx,T^ny)<\epsilon$ for all $n \in \Z_+$. A transitive
TDS is called {\em almost equicontinuous} if it has at least one
equicontinuity point. If a transitive system
is almost equicontinuous then the set of equicontinuity points
coincides with the set of transitive points and hence it is dense
$G_{\delta}$. A transitive TDS $(X,T)$ is called {\em sensitive}
if there is an $\ep>0$ such that whenever $U$ is a nonempty open set there
exist $x,y\in U$ such that $d(T^nx,T^ny)>\ep$ for some $n\in \N$.
A transitive TDS is either almost equicontinuous or sensitive.
In particular a minimal system is either equicontinuous or sensitive
(see \cite{GW93} and \cite{AAB}).

A TDS $(X,T)$ is said to be {\em chaotic in the sense of Devaney}
(or an infinite {\em $P$-system}) if it is transitive and $X$ is infinite with a dense set of periodic points. Such a system is always sensitive
(see \cite{BBC} and \cite{GW93}).

\medskip

\subsection{Families and filters}
We say that a collection $\F$ of subsets of $\Z_+$ (or $\Z$) is
a {\em a family} if it is hereditary upward, i.e. $F_1 \subseteq F_2$
and $F_1 \in \F$ imply $F_2 \in \F$. A family $\F$ is called {\em proper} if it is neither empty nor the entire power set of $\Z_+$, or, equivalently if 
$\Z_+ \in \F$ and $\emptyset \not\in \F$.
If a family $\F$ is closed under finite intersections and is proper,
then it is called a {\em filter}. 
A collection of nonempty subsets $\mathcal{B}$ is a {\em filter base}
if for every $B_1,B_2 \in \mathcal{B}$ there is $B_3 \in \mathcal{B}$
with $B_3 \subset B_1 \cap B_2$. When $\mathcal{B}$
is a filter base then the family 
$$
\mathcal{F}=\{F : \exists B \in \mathcal{B}\ \text{with}\ B \subset F\},
$$
is a filter. A maximal filter is called an {\em ultrafilter}.
By Zorn's lemma every filter is contained in an ultrafilter.

For a family $\F$ its {\em dual} is the family
$\F^{\ast} :=\{F\subseteq \Z_+ | F \cap F' \neq \emptyset \ for \
all \ F' \in \F \}$. Any nonempty collection $\A$ of subsets  of $\Z_+$
generates a family $\F(\A) :=\{F \subseteq \Z_+:F \supset A$ for some 
$A \in \A\}$.

The collection $\beta\Z$ of ultrafilters on $\Z$ can be identified 
with the \v{C}ech-Stone compactification of the integers,
where to $n \in \Z$ corresponds the principle unltrafilter
$\{A: n \in A \subset \Z\}$. Using the universal property of
this compactification one shows that the map $n \mapsto n+1$
on $\Z$ extends to a homeomorphism $S: \beta\Z \to \beta\Z$
and that, more generally, addition in $\Z$ can be extended to
a binary operation on $\beta\Z$ making it an Ellis semigroup; 
i.e.  for every $p \in \beta\Z$, right multiplication
$R_p: q \mapsto qp$ is continuous. In fact the resulting dynamical system
$(\beta\Z,S)$ is the universal point transitive dynamical system
and the corresponding enveloping semigroup is naturally
identified with $\beta\Z$ itself via the map $p \mapsto L_p$,
where $L_p: q \mapsto pq$. A similar construction defines
$\beta\Z_+$.
In view of these facts the Ellis semigroup $\beta\Z$ (or $\beta\Z_+$)
can serve as a {\em universal enveloping semigroup}
(see Subsection \ref{es} above).

\medskip

\begin{lem}\label{lem-tran}
Let $(X,T)$ be a transitive TDS. Then the collection of sets
$$
\A = \{N(U,U): \text{$U$ is a nonempty open subset of $X$} \}
$$
is a filter base, whence the family $\F(\A)$ is a filter.
\end{lem}

\begin{proof}
Let $U_1$ and $U_2$ be nonempty open subsets of
$X$. As $(X,T)$ is transitive, there is an $n\in \N$ such that $U_3=U_1\cap
T^{-n}U_2\neq \emptyset$. Then
\begin{eqnarray*}
N(U_3,U_3)& \subseteq & N(U_1,U_1)\cap N(T^{-n}U_2,T^{-n}U_2)\\
&= &N(U_1,U_1)\cap N(T^{n}T^{-n}U_2,U_2)\\ &\subseteq &
N(U_1,U_1)\cap N(U_2,U_2),
\end{eqnarray*}
and our claim follows.
\end{proof}

For a TDS $(X,T)$ and a point $x \in X$ define 
$$
\mathcal{I}_x = \{ N(x,U): 
\text{$U$ is a neighborhood of $x$}\}.
$$ 
A point
$x$ is {\em recurrent} for $(X,T)$ if and only if  each such return time set
$N(x,U)$ is nonempty and so if and only if  $\mathcal{I}_x$  is a filter base. 
For a pair $(x_1,x_2) \in X \times X$ define 
$$
\P_{(x_1,x_2)} = \{ N((x_1,x_2),V) : 
\text{$V$ is a neighborhood of the
diagonal in $X \times X$}\}.
$$ 
A pair $(x_1,x_2)$ is proximal if and only if  each such $N((x_1,x_2),V)$ is nonempty and so
if and only if  $\P_{(x_1,x_2)}$ is a filter base.

\medskip

\section{Transitivity, rigidity and proximality}\label{Sec-tran}

\subsection{Rigid and proximal sets}
The following definitions are from \cite{GM}, where they were
defined for the total space $X$.

\begin{de}\label{rigid}
Let $(X,T)$ be a TDS , $K\subseteq X$ and $S \subset \Z_+$.
\begin{enumerate}
\item
We say that $K$ is {\em rigid} with respect to a sequence
$S =\{n_k\}_{k=1}^\infty$, $n_k\nearrow +\infty$ if 
$\lim \limits_{k\to \infty}T^{n_k}x=x$ for every $x\in K$.
\item
$K$ is {\em uniformly rigid} with respect to $S$ if for every $\ep>0$
there is an $n \in S$ with $d(T^n x, x) < \ep$ for all $x$ in $K$.
\item
$K$ is {\em weakly rigid} with respect to $S$ if
every finite subset of $ K$ is uniformly rigid with respect to $S$.
\end{enumerate}
\end{de}

In items (2) and (3) we omit the reference to $S$
when $S = \Z_+$. Clearly 
$$
\text{\rm{uniform rigidity $\Rightarrow$ rigidity $\Rightarrow$
weak rigidity.}}
$$

\medskip

Recall that the $A(X,T)$, the adherence semigroup of $(X,T)$,
is defined as the $\omega$-limit set of the
collection $\{T^n: n \in \Z_+\}$ in $E(X,T)$.
We have the following lemma (see \cite{GM}).

\begin{lem}
A subset  $K \subset X$ is weakly rigid if and only if there is an 
idempotent $u \in A(X,T)$ with $ux =x$ for every $x \in K$. 
In particular the identity map ${\rm id}: X \to X$ is an element
of $A(X,T)$ if and only if  the system $(X,T)$ is weakly rigid.
\end{lem}

\begin{proof}
A subset $K\subset X$ is weakly rigid if and only if  
$\bigcup \{ \mathcal{I}_x  : x \in K \}$ is a filter base and so is contained in
some ultrafilter. This implies that when $K$ is weakly rigid the set
$$
S_K = \{p \in A(X,T): px=x\ \text{for every } x \in K\}
$$
is a nonempty closed subsemigroup of $A(X,T)$. By Ellis' lemma
there is an idempotent $u \in S_K$.
The converse is clear.
\end{proof}


%


%


\begin{rem}
We let for $n \ge 1$,
$$
Recur_n(X) = \{(x_1,\dots,x_n) \in X^n: \forall \ep >0,\ \exists
k \in \Z_+ \ {\text{with}}\ d(T^k x_i, x_i ) < \ep, \forall i\}.
$$
In this notation $K \subset X$ is weakly rigid if and only if
for every $n$, every $n$-tuple $(x_1,\dots,x_n) \in K^n$
is in $Recur_n(X)$. 
Note that $Recur_n(X)$ is a $G_\delta$ subset of $X^n$.
\end{rem}

\begin{de}\label{proximal}
Let $(X,T)$ be a TDS, $K \subset X$ and $S \subset \Z_+$.
\begin{enumerate}
\item
A subset $K$ of $X$ is called {\em pairwise proximal}
if every pair $(x,x') \in K \times K$ is proximal.
\item
The subset $K$ is called {\em uniformly  proximal} with respect to $S$ 
if for every $\ep>0$ there is $n \in S$ with ${\rm diam}\, T^n K < \ep$.
\item
A subset $K$ of $X$ is called {\em proximal} with
respect to  $S$ if every finite subset of $K$ is
{\em uniformly } proximal with respect to $S$.

\end{enumerate}
\end{de}

\begin{rem}
Thus, $K$ is uniformly  proximal with respect to $S$ when there is a sequence
$\{ n_k \}$ in $S$ such that ${\rm diam}\, T^{n_k} K $ converges to $0$.
A subset $K\subset X$ is proximal if and only if 
$\bigcup \{ \P_{(x_1,x_2)}  : (x_1,x_2) \in K \times K \}$
is a filter base and so is contained in
some ultrafilter.  It follows that $K \subset X$ is
a proximal set if and only if there exists an element $p \in E(X,T)$
with $pK = \{x\}$ for some $x\in X$.
We let for $n \ge 1$, 
\begin{align*}
Prox_n(X) = & \{ (x_1,x_2,\cdots,x_n): \forall  \epsilon>0 \ \exists
m\in \N  \text{ such that }\\
& \hskip0.3cm  \text{diam}\,(\{ T^mx_1,\cdots,T^mx_n \})<\epsilon \}.
\end{align*}
In this notation $K \subset X$ is a proximal set if and only if
for every $n$, every $n$-tuple $(x_1,\dots,x_n) \in K^n$
is in $Prox_n(X)$. 
Again we note that $Prox_n(X)$
is a $G_\delta$ subset of $X^n$.
\end{rem}



\medskip

\subsection{Transitivity implies partial rigidity}


A nonempty subset $K$ of a compact space $X$ is a {\em Mycielski} set
if it is a countable union of Cantor sets. In the following theorem we show
that every transitive TDS contains a dense weakly rigid Mycielski subset.
While we will later derive this result, and more, from the Kuratowski-Mycielski Theorem, we include here a direct proof which employs an explicit
construction rather than an abstract machinery (see Theorem \ref{tran'}
below).

\begin{thm}\label{tran}
Let $(X,T)$ be a transitive TDS without isolated points. Then
there are Cantor sets $C_1 \subseteq C_2 \subseteq \cdots$ such
that $K=\bigcup \limits_{i=1}^{\infty}C_n$ is a dense rigid
subset of $X$ and for each $N \in \N$, $C_N$ is uniformly rigid.

If in additional, for each $n \in \N$, $Prox_n(X)$ is dense in
$X^n$, then we can require that for each $N \in \N$, $ C_N$ is
uniformly  proximal, whence $K$ is a proximal set.
\end{thm}

\begin{proof} Let $Y= \{y_1,y_2,\dots \}$ be a countable dense
subset of $X$ and for each $n \ge 1$ let
$Y_n= \{y_1,y_2, \dots, y_n\}$.
Let $\F$ be the smallest family containing the collection
$$
\{N(U,U): \text{$U$ is a nonempty open subset of $X$} \}.
$$
Since $(X,T)$ is transitive, $\F$ is a filter by Lemma \ref{lem-tran}.
Let $a_0=0$ and $V_{0,1}=X$. We have the following claim.

\medskip

\noindent{\bf Claim:}
For each $S\in \F^*$ there are sequences
$\{a_n\}\subseteq \N$, $\{k_n\} \subseteq S$, and sequences $\{U_n\}_{n=1}^\infty$ and
$\{V_{n,1},V_{n,2},\cdots, V_{n,a_n}\}_{n=1}^{\infty}$ of nonempty open subsets of $X$ with the following properties:
\begin{enumerate}
\item $2a_{n-1} \le a_n \le 2a_{n-1}+n$.
\item $diam\, V_{n,i} < \frac 1n,$ $ i=1,2,\dots,a_n$.
\item The closures $\{ \ol {V_{n,i} } \}_{i=1}^{a_n}$ are pairwise
disjoint.
\item $\ol {V_{n,2i-1}} \cup \ol {V_{n,2i}} \subset V_{n-1,i}$,
$i=1,2,\cdots,a_{n-1}$.
\item $Y_n \subset B(\bigcup \limits_{i=1}^{a_n} V_{n,i},\frac 1n)$,
where $B(A,\ep):=\{x \in X: d(x,A)<\ep \}$ .
\item $T^{k_n}(V_{n,2i-1}\cup V_{n,2i}) \subseteq V_{n-1,i}$,
$i=1,2, \cdots, a_{n-1}$.
\end{enumerate}

\medskip

\noindent{\it Proof of Claim:}
For $j=1$, take $a_1 = 2,\  k_1 =1$, and $V_{1,1},V_{1,2}$
any two nonempty open sets of diameter $< 1$ with disjoint closures
such that $y_1 \in B(V_{1,1} \cup V_{1,2}, 1)$.
Suppose now that for $1\le j \le n-1$ we have $\{a_j\}_{j=1}^{n-1}$,
$\{k_j\}_{j=1}^{n-1}$ and $\{V_{j,1},V_{j,2},\cdots,
V_{j,a_j}\}$, satisfying conditions (1)-(6).

Choose $2a_{n-1} \le a_n \le 2a_{n-1}+n$ and nonempty open subsets
$V^{(0)}_{n,1}, V^{(0)}_{n,2}, \cdots, V^{(0)}_{n,a_n}$ of $X$
such that:

(a) $diam\, V^{(0)}_{n,i}< \frac 1{2n}$, $i=1,2,\cdots,a_n$.

(b) The closures $\{ \ol {V^{(0)}_{n,i}}\}_{i=1}^{a_n}$ are
pairwise disjoint.

(c) $\ol {V^{(0)}_{n,2i-1}} \cup \ol {V^{(0)}_{n,2i}} \subset
V_{n-1,i},\ i=1,2,\cdots,a_{n-1}$.

(d) $Y_n \subset B(\bigcup \limits_{i=1}^{a_n} V^{(0)}_{n,i},\frac
1{2n})$.

As $N(V^{(0)}_{n,i},V^{(0)}_{n,i})\in \F$ for each $1\le i\le
a_n$, $\bigcap \limits_{i=1}^{a_n} N(V^{(0)}_{n,i},V^{(0)}_{n,i})
\in \F$. Take $k_n \in S \cap \bigcap \limits
_{i=1}^{a_n}N(V^{(0)}_{n,i},V^{(0)}_{n,i})$. Hence there are
nonempty open sets $V^{(1)}_{n,i} \subseteq V^{(0)}_{n,i}$, $1\le i \le a_n$,
such that

(e) $T^{k_n}(V^{(1)}_{n,2i-1}\cup V^{(1)}_{n,2i}) \subseteq
V_{n-1,i}$, $i=1,2, \cdots, a_{n-1}$.

Let $V_{n,i}=V^{(1)}_{n,i}$, $1\le i \le a_n$. Then  the
conditions $(1)-(6)$ hold for $n$. By induction we have the
claim.

\medskip

Let $C_n= \bigcap \limits_{j=n}^{\infty} \bigcup
_{i=1}^{2^{j-n}a_n} \ol{V_{j,i}}$. Then $C_1 \subseteq C_2
\subseteq \cdots$, and by $(1)-(4)$, $C_n$ is a Cantor set.
By (2),(4) and (5), $K=\bigcup \limits_{n=1}^{\infty} C_n$ is dense in $X$.
For each $N\in \N$, by (6),
$C_N$ is uniformly rigid with respect to a subsequence of $S$.

Finally, if in addition for each $n\in \N$, $Prox_n(X)$ is dense in
$X^{(n)}$ then we can in the above construction, when choosing the subsets
$V^{(1)}_{n,i} \subseteq V^{(0)}_{n,i}$, $1\le i \le a_n$,
add the following condition to the claim above:
$$
(7)\  \text{for each $n\in \N$ there is $t_n\in \N$ such that}\ {\rm diam}\,
T^{t_n}(\bigcup_{i=1}^{a_n}
\overline{V_{n,i}})<\frac{1}{n}.
$$
By the requirement $(7)$ we obtain that for each $N\in \N$,
$C_N$ is uniformly  proximal with respect to $\{ t_n \}$.
\end{proof}

\medskip
\begin{rem}
Using the fact that the set $X_{tr}$ of transitive points is dense in $X$
we can in the above construction, when choosing the subsets
$V^{(1)}_{n,i} \subseteq V^{(0)}_{n,i}$, $1\le i \le a_n$,
add the following condition to the claim in the proof:
$$
(8)\ Y_n \subseteq B(Orb(x,T),\frac 1n)\ \text{ for each}\ x\in \bigcup
\limits_{i=1}^{a_n}\overline {V_{n,i}}.
$$
It then follows that every point in $\bigcup \limits_{i=1}^{\infty}C_n$ is a
transitive point.
\end{rem}

\medskip

Motivated by Theorem \ref{tran} we define uniformly  chaotic set as
follows:
\begin{de}
Let $(X,T)$ be a TDS.
A  subset $K \subseteq X$ is called a {\em uniformly  chaotic set}
if there are Cantor sets $C_1 \subseteq C_2 \subseteq \cdots$ such that
\begin{enumerate}
\item
$K = \bigcup \limits_{i=1}^{\infty}C_n$ is a rigid subset of
$X$ and also a proximal subset of $X$;
\item
for each $N\in \N$, $C_N$ is
uniformly rigid; and
\item
for each $N\in \N$, $C_N$ is uniformly  proximal.
\end{enumerate}
$(X,T)$ is called {\em (densely) uniformly chaotic}, if $(X,T)$ has
a (dense) uniformly  chaotic subset.
\end{de}


\begin{rem}
Actually the fact that $K$ is rigid and proximal follows from the 
conditions (2) and (3):
Let $J_N  = \{ n :  d(T^nx,x) < 1/N \ \text{for all}\  x \in \bigcup_{i=1}^N \ C_i\}$.
Each $J_N$ is nonempty by assumption and
$J_{N+1} \subset J_N$.  Choose $n_N \in J_N$.  As $N \to \infty$,
$T^{n_N}x \to x$ for all $x \in \bigcup_{i=1}^{\infty} \ C_i$.
Thus $K$ is a rigid subset with respect to the sequence $\{n_N\}$.
Clearly condition (3) implies that $K$ is a proximal set. 
\end{rem}

\medskip

Obviously, a uniformly chaotic set is an uncountable strongly scrambled
set, hence every uniformly chaotic system is strongly Li-Yorke chaotic.
Restating Theorem \ref{tran} we have:
\begin{thm} \label{unc}
Let $(X,T)$ be a transitive TDS without isolated points.
If for each $n\in \N$, $Prox_n(X)$ is dense in $X^{(n)}$, then
$(X,T)$ is densely uniformly chaotic. In particular every such system
is strongly Li-York chaotic.
\end{thm}

\medskip

\section{A criterion for chaos and applications}\label{Sec-cri}

\subsection{A criterion for chaos}

\begin{thm}[A criterion for chaos] \label{zhunzhe}
Let $(X,T)$ be a transitive TDS without isolated points. If there
is some subsystem $(Y,T)$ of $(X,T)$ such that $(X\times Y,T)$ is
transitive, then $(X,T)$ is densely uniformly chaotic.
\end{thm}

\begin{proof}
By Theorem \ref{unc}, it suffices to show that
for each $n\in \N$, $Prox_n(X)$ is dense in $X^{(n)}$. For a fixed
$n\in \N$ and any $\ep >0$ let
$$
P_n({\ep})=\{(x_1,x_2,\dots,x_n): \exists
m\in \N  \text{ such that }\\
\text{diam}\,(\{ T^mx_1,\dots,T^mx_n \})<\epsilon \}.
$$
Thus $Prox_n(X)=\bigcap_{m=1}^{\infty}P_n(\frac{1}{m})$
and by Baire's category theorem it is enough to show that
for every $\ep>0$, $P_n({\ep})$ is a dense open subset of $X^n$.

Fix $\ep>0$, let $U_1,U_2,\cdots,U_n$ be a sequence of
nonempty open subsets of $X$, and let $W$ be a nonempty open
subset of $Y$ with $\text{diam}\,(W)<\ep$.
By assumption $(X\times Y,T)$ is transitive, whence
$$
N(U_1 \times W,U_2\times W) = N(U_1,U_2)\cap N(W\cap Y,W\cap Y)\neq \emptyset.
$$
Let $m_2$ be a member of this intersection. Then
$$
U_1\cap T^{-m_2}U_2\neq \emptyset\quad
 \text{and}\quad W\cap T^{-n_2}W \cap Y\neq \emptyset.
$$
By induction, we choose natural numbers $m_3,m_4,\cdots,m_n$
such that
$$
U_1\cap \bigcap_{i=2}^n T^{-m_i}U_i\neq \emptyset\quad  \text{and}\quad
W\cap \bigcap_{i=2}^n T^{-m_i}W \cap Y\neq \emptyset.
$$

Since $(X,T)$ is transitive, there is a transitive point
$x\in U_1\cap \bigcap_{i=2}^n T^{-m_i}U_i$ and let $y \in W\cap
\bigcap_{i=2}^n T^{-m_i}W$. Since $x$ is a transitive point, there
exists a sequence $l_k$ such that  $\lim_{k\to \infty}
T^{l_k}x=y$. Thus, $\lim_{k\to \infty} T^{l_k}(T^{m_i}x)=T^{m_i}y$ for each
$2\le i\le n$. Since $\{ y,T^{m_2}y,\dots,T^{m_3}y \}\subset W$
and $\text{diam}\,(W)<\ep$, for large enough $l_k$, we have
$$
\text{diam}\,(\{ T^{l_k}x,T^{l_k}(T^{m_2}x),\dots,T^{l_k}(T^{m_n}x)\})<\ep.
$$
That is, $(x,T^{m_2}x,\dots,T^{m_n}x)\in P_n(\ep)$. Noting that
$(x,T^{m_2}x,\dots,T^{m_n}x)\in U_1\times U_2\times \cdots \times
U_n$, we have shown that
$$
P_n(\ep)\cap U_1\times U_2\times \cdots \times U_n \neq
\emptyset.
$$
As $U_1,U_2,\cdots,U_n$ are arbitrary, $P_n({\ep})$ is indeed
dense in $X^{(n)}$.
\end{proof}

\medskip

\subsection{Some applications}

In the rest of this section we will obtain some applications of the above
criterion. First, we need to recall some definitions (see
\cite{BHM,HY02}).

\medskip

Two topological dynamical systems are said to be {\em weakly disjoint} if their product is transitive. Call a TDS $(X,T)$:
\begin{itemize}
\item
{\em scattering} if it is weakly disjoint from every minimal system;
\item
{\em weakly scattering} if it is weakly disjoint from every minimal equicontinuous system;
\item
{\em totally transitive} if it is weakly disjoint from every periodic system.
(Check that this is equivalent to the usual definition which
requires that $(X,T^n)$ be transitive for all $n \ge 1$.)
\end{itemize}

\medskip

Using this terminology and applying Theorem \ref{unc}
we easily obtain the following:
\begin{cor} \label{app}
If $(X,T)$ is a TDS without isolated points and one of the following properties,
then it is densely uniformly  chaotic:
\begin{enumerate}
\item $(X,T)$ is transitive and has a fixed point;
\item $(X,T)$ is totally transitive with a periodic point;
\item $(X,T)$ is scattering;
\item $(X,T)$ is weakly scattering with an equicontinuous minimal subset;
\item $(X,T)$ is weakly mixing.
\end{enumerate}
Finally
\begin{enumerate}
\item[(6)]
If $(X,T)$ is transitive and has a periodic point of order $d$,
then there is a closed $T^d$-invariant subset $X_0 \subset X$,
such that $(X_0,T^d)$ is densely uniformly  chaotic and
$X = \bigcup_{j=0}^{d-1} T^j X_0$. In particular $(X,T)$
is uniformly  chaotic.
\end{enumerate}
\end{cor}

\begin{proof}
The only claim that needs a proof is (6).
Suppose $y_0\in X$ is a periodic point of period $d$ and
let $x_0$ be a transitive point; so that $\overline{Orb_T(x_0)} = X$.
Set $\overline{Orb_{T^d}(x_0)} = X_0$ (this may or may not be all
of $X$). In any case the dynamical system $(X_0, T^d)$ is transitive,
and has a fixed point. Thus, by case (1), it is densely uniformly  chaotic
for $T^d$. Both uniform proximality
and uniform rigidity of subsets go over to $(X,T)$, hence
$(X,T)$ is uniformly  chaotic.
Clearly $X = \bigcup_{j=0}^{d-1} T^j X_0$.
\end{proof}

Part (6) provides a new proof of a result of J-H. Mai \cite{Mai},
and as in Mai's paper we have the following corollary.

\begin{thm}
Devaney chaos implies uniform chaos.
\end{thm}

One can strengthen Corollary \ref{app} (2) in the following way.

\begin{de}
Let $(X,T)$ be a TDS. A point $x\in X$
is  {\em regularly almost periodic} if for each neighborhood $U$ of $x$
there is some $k\in\N$ such that $k\Z_+ \subseteq N(x,U)$.
Note that such a point is in particular a minimal point
(i.e. its orbit closure is minimal).
\end{de}

\begin{rem}
Let $(X,T)$ be a minimal system. Then $(X,T)$
contains a regularly almost periodic point if and only if it is an
almost one-to-one extension
of an adding machine. If in addition $(X,T)$ is a subshift then it is
isomorphic to a Toeplitz system (see e.g. \cite{MP}).
\end{rem}

Next we recall the following definition from \cite{AG}.

\begin{de}
A property of topological dynamical systems is said
to be {\em residual} if it is non-vacuous and is inherited by factors,
almost one-to-one lifts, and inverse limits.
\end{de}

It is not hard to check that being weakly disjoint from a fixed TDS
$(X,T)$ is a residual property (see \cite{AG}).
One can also show that the smallest class of TDS which contains the periodic
orbits and is closed under inverse limits and almost one-to-one
extensions is exactly the class of almost one-to-one extensions of
adding machines.
It now follows that a TDS is totally transitive if and only if 
it is weakly disjoint from every almost one-to-one extension of an
adding machine.

\begin{cor}
If $(X,T)$ is totally transitive with a regularly almost periodic point,
then it is densely uniformly  chaotic.
\end{cor}

In a similar way we see that a TDS is weakly scattering
if and only if  it is weakly disjoint from every system which is an almost one-to-one
extension of a minimal equicontinuous system (these systems
are also called {\em almost automorphic}). Thus we also have
a stronger version of Corollary \ref{app} (4)

\begin{cor}
If $(X,T)$ is weakly scattering and has an almost automorphic subsystem
then it is densely uniformly  chaotic.
\end{cor}

The following example shows that we can not weaken the condition
``total transitivity" to ``transitivity".

\medskip

\noindent {\bf Example.}
Let $(X,T)$ be a Toeplitz system and let $\pi : X \to Z$ be the corresponding
almost one-to-one factor map from $X$ onto its maximal adding machine
factor. Clearly then every proximal set in $X$ is contained in a fiber
$\pi^{-1}(z)$ for some $z \in Z$. Suppose now that $|\pi^{-1}(z)| < \infty$
for every $z \in Z$, and that for some $z \in Z$ there are points
$x, y \in \pi^{-1}(z)$ such that $(x,y)$ is a recurrent pair and
therefore a strong Li-Yorke pair
(one can easily construct such systems, see e.g. \cite{W}).
Let $Y=\overline {Orb((x,y),T)} \subseteq X\times X$.
By assumption the point $(x,y)$ is recurrent in $X \times X$
and forms a proximal pair.
Thus the system $(Y,T)$ is transitive and, as one can easily check, has
$\Delta_{X}$ as its unique minimal subset. Since $(X,T)$ is an almost
one to one extension of an adding machine the diagonal
$\Delta_{X}\subset Y$ contains regularly almost periodic points.
However $(Y,T)$ can not be Li-Yorke chaotic because
our assumption implies that every proximal set in $Y$ is finite.

This example also shows the existence of a non-minimal transitive system
which is not Li-Yorke chaotic.

\medskip

\section{The Kuratowski-Mycielski Theory}\label{Sec-KM}

Let $X$ be a compact metric space. We recall that a subset $A\subset X$ is
called a {\em Mycielski set} if it is a union of countably many Cantor sets.
(This definition was introduced in \cite{BGKM}. Note that in \cite{Ak03}
a Mycielski set is required to be dense.)
The notion of independent sets and the corresponding topological machinery
were introduced by Marczewski \cite{Mar}, and Mycielski \cite{M}.
This theory was further developed by Kuratowski in \cite{K}. 
The first application to dynamics is due to Iwanik \cite{I}. 
Consequently it was used as a main tool in
\cite{BGKM}, where among other results the authors showed that
positive entropy implies Li-Yorke chaos.
See \cite{Ak03} for a comprehensive treatment of this topic.

In this section we first review the Kuratowski-Mycielski theory,
mainly as developed in \cite{Ak03}, and then consider the results
of Sections \ref{Sec-tran} and \ref{Sec-cri} in view of this theory.

\medskip

\subsection{The Kuratowski-Mycielski Theorem}

We begin by citing two classical results.

\begin{thm}[Ulam]\label{u}
Let $\phi: X\to Y$ be a continuous open surjective map with $X$ and $Y$
metric compact spaces. If $R$ is a dense $G_{\delta}$ subset of
$X$, then
\begin{equation*}
Y_0=\{y \in Y: \phi^{-1}(y)\cap R \ \text{is dense in } \
\phi^{-1}(y)\}
\end{equation*}
is a dense $G_{\delta}$ subset of $Y$.
\end{thm}

\begin{thm}[Mycielski]\label{m}
Let $X$ be a complete metric space with no isolated points. Let
$r_n \nearrow \infty$ be a sequence of positive integers and for every $n$ let
$R_n$ be a meager subset of $X^{r_n}$. Let
${\{O_i\}}_{i=1}^{\infty}$ be a sequence of nonempty open subsets of $X$.
Then there exists a sequence of Cantor sets $C_i \subset O_i$ such that
the corresponding  Mycielski set $K={\bigcup}_{i=1}^{\infty} C_i$ has the
property that for every $n$ and every $x_1,x_2, \dots, x_{r_n}$, distinct elements of
$K$, $(x_1,x_2,\dots,x_{r_n})\not \in R_n$.
\end{thm}

An especially useful instance of Mycielski's theorem is obtained
as follows (see \cite{Ak03}, Theorem 5.10, and \cite{AAG}, Theorem 6.32).
Let $W$ be a symmetric dense $G_\delta$ subset
of $X \times X$ containing the diagonal $\Delta_X$, and let
$R = X \times X \setminus W$.
Let $r_n=n$ and set
$$
R_n =\{(x_1,\dots,x_n): (x_i, x_j) \not\in W,\ \forall\  i \neq j\}.
$$

\begin{thm}\label{special}
Let $X$ be a perfect compact metric space and $W$
a symmetric dense $G_\delta$ subset
of $X \times X$ containing the diagonal $\Delta_X$.
There exists a dense Mycielski subset $K \subset X$ such that $K \times
K \subset W$.
\end{thm}

\medskip

We collect some notation and results from Akin \cite{Ak03}.
For $X$ a compact metric space we denote by $C(X)$ the compact
space of  closed subsets of $X$ equipped with the Hausdorff metric.
Since $\emptyset$ is an isolated point,
$C'(X) = C(X) \setminus \{ \emptyset \}$ is compact as well.

We call a collection of sets $Q \subset C'(X)$ {\em hereditary}
if it is hereditary downwards, that is,
$A \in Q$ implies $C'(A) \subset Q$ and, in particular,
every finite subset of $A$ is in $Q$.
For a hereditary subset $Q$ we define
$R_n(Q) = \{ (x_1, \dots ,x_n) \in X^n :
\{ x_1, \dots ,x_n \} \in Q \} = i_n^{-1}(Q)$ where $i_n : X^n \to C'(X)$
is the continuous map defined by
$i_n(x_1, \dots ,x_n) = \{ x_1, \dots ,x_n \}$.
In particular, if $Q$ is a $G_{\delta}$ subset of
$C'(X)$ then $R_n(Q)$ is a $G_{\delta}$ subset of $X^n$ for all $n$.
Call $A$ a $\{ R_n(Q) \}$ set if $A^n \subset R_n(Q)$ for
all $n = 1,2, \dots $ or, equivalently, if every finite subset of $A$ lies in $Q$.
Clearly, the union of any chain of $\{ R_n(Q) \}$ sets is an $\{ R_n(Q) \}$
set and so every $\{ R_n(Q) \}$ set is contained in a maximal $\{ R_n(Q) \}$ set.

If $D \subset X$ we define
$Q(D) = \{ A \in C'(X) : A \subset D \}$, for which $R_n = D^n$. If
$B \subset X \times X$ is a subset which satisfies
$$
(x,y) \in B \qquad \Longrightarrow \qquad (y,x), (x,x) \in B
$$
then we define
$Q(B) = \{ A \in C'(X) : A \times A \subset B \}$, for which $R_2 = B$ and
$(x_1, \dots ,x_n) \in R_n$ if and only if  $(x_i,x_j) \in B$ for all $i,j = 1, \dots ,n$.

If $D$ (or $B$) is a $G_{\delta}$ then so is $Q(D)$ (resp. $Q(B)$).
Because the finite sets are dense in $C'(X)$ it follows that if $D$ is  dense
in $X$ (or $B$ is $G_{\delta}$ and dense in $X \times X$) then $Q(D)$
(resp. $Q(B)$ ) is dense in $C'(X)$.

\medskip

{\bfseries Examples:}
(1) Let $Q(Recur) = \{ A \in C'(X) : A\ \text{is uniformly rigid}\}$.
We denote by $Recur_n$ the set
$R_n(Q(Recur)) = \{ (x_1, \dots ,x_n) : \text{recurrent in}\ X^n \}$.
The $\{ Recur_n \}$ subsets are the weakly rigid subsets.
For fixed $n$ and $ \ep$ the condition $d(T^n x, x) < \ep$ for  all $x \in A$ is
an open condition on $A \in C'(X)$.  Hence, $Q(Recur)$ and $Recur_n$ are $G_{\delta}$ sets.

Notice that if $x$ is a transitive point for a transitive TDS $(X,T)$ then
points of the form $(T^{k_1}x, \dots ,T^{k_n}x)$ comprise a dense set of
recurrent points in $X^n$.  Thus, for a transitive system $Recur_n$ is
dense in $X^n$.
In addition $X_{tr}$ is a dense $G_{\delta}$ in $X$ and so
$Q(X_{tr}) = \{ A \in C'(X) : A \subset X_{tr} \}$ is a dense $G_{\delta}$
subset of $C'(X)$.

\medskip

(2) Let $Q(Prox) = \{ A \in C'(X) : A\ \text{is uniformly  proximal}\}$.
We denote by $Prox_n$ the set
$R_n(Q(Prox))$.  The $\{ Prox_n \}$ subsets are the proximal subsets.
For fixed $n$ and $\ep$ the condition ${\rm diam}\, T^n A < \ep$ is
an open condition on $A \in C'(X)$.  Hence, $Q(Prox)$ and $Prox_n$ are $G_{\delta}$ sets.

$Prox_2 = P(X,T)$ the set of proximal pairs.
The $G_{\delta}$ set $Q(P(X,T))$ is the set of compacta $A$ such
that $A \times A \subset P(X,T)$. The $\{ R_n(Q(P)) \}$ sets are the pairwise proximal sets.

\medskip

(3) For use below we define for $Y$ a closed subset of $X$:
\begin{align*}
Q(TRANS,Y) = \{& A \in C'(X):\  \text{for every}\ \ep > 0, \ n \in \Z_+,
\ \text{pairwise disjoint closed}\\
& A_1, \dots ,A_n \subset A \
\text{and} \ y_1, \dots ,y_n \in Y,\ \text{there exists a positive}\\
& \text{integer $k$ such that}\ d(T^k x,y_i) < \ep\
\text{for all}\ x\in A_i,\ i = 1, \dots ,n \}.
\end{align*}
It is easy to check that $Q(TRANS,Y)$ is a $G_{\delta}$ set,
see Akin \cite{Ak03}, Lemma 6.6(a). Clearly,
$(x_1, \dots ,x_n) \in R_n(Q(TRANS,Y))$ if and only if  for every $\ep > 0$ and
$y_1, \dots ,y_n \in Y$ there exists $k$
such that $d(T^k x_i, y_i) < \ep$ for $i = 1, \dots ,n$.

\medskip

The point of the peculiar condition is given by

\begin{lem}\label{lem-transY}
If $K$ is a Cantor set in $X$, then $K \in Q(TRANS,Y)$ if and only if  for every
continuous map $h : K \to Y$ and every $\ep > 0$ there exists a positive
integer $k$ such that $d(T^k x,h(x)) < \ep$ for all $x \in K$.
\end{lem}

\begin{proof}
Recall that the locally constant functions on $K$, which are the continuous functions with finite range, form a
dense subset of $\mathcal{C}(K,Y)$ the space of continuous functions.
It thus suffices to consider such functions $h$.
If $h(K)$ is the set $ \{ y_1, \dots ,y_n \}$ of $n$ distinct points then
$\{ A_i = h^{-1}(y_i) : i = 1, \dots ,n \}$ is
a clopen partition of $K$. Hence, $K \in Q$ implies there exists a $k$ such that $T^k \rest K$ is within $\ep$ of $h$.

Conversely, given disjoint closed sets $A_1, \dots ,A_n$ in
$K$ and points $y_1, \dots ,y_n \in Y$ there exists a
clopen partition $B_1, \dots ,B_n$ of $X$ with $A_i \subset B_i$ for
$i = 1, \dots ,n$.
The function $h : X \to Y$ with $h(x) = y_i$ for $x \in B_i$ is continuous and approximating it by some $T^k \rest K$ shows that $K \in Q$.
\end{proof}

For a TDS $(X,T)$ and closed $Y \subset X$, motivated by Lemma
\ref{lem-transY}, we will call a Cantor set $K \in Q(TRANS,Y)$ a
{\em Kronecker set for $Y$}. 

\begin{lem}\label{Kronecker}
Let $(X,T)$ be a TDS and $Y$ a closed nonempty subset of $X$.
Then any Kronecker set for $Y$,  $K \in Q(TRANS,Y)$, is uniformly proximal.
If moreover $K \subset Y$ then $K$ is also uniformly rigid, 
hence uniformly chaotic.
\end{lem}

\begin{proof}
Apply Lemma \ref{lem-transY}.
For the first assertion take $h: K \to Y$ as any constant map
$h: K \to Y$, $h(x) = y_0, \ \forall  x \in K$. For the second, take
$h: K \to Y$ as $h(x) = x, \ \forall x \in K$. 
\end{proof}

If $X$ is a perfect, nonempty, compact metric space then $CANTOR(X)$
the set of Cantor sets in $X$ is a dense $G_{\delta}$ subset of
$C'(X)$, see e.g.  Akin \cite{Ak03} Propsition 4.3(f).

\medskip

The importance of all this stems from the {\em Kuratowski-Mycielski Theorem}.  This version comes from Akin \cite{Ak03} Theorem 5.10 and Corollary 5.11.

\begin{thm}\label{KMT}
For $X$ a perfect, nonempty, compact metric space, let $Q$ be a $G_{\delta}$ subset of $C'(X)$.

(a) The following conditions are equivalent
\begin{enumerate}
\item For $n = 1,2, \dots , \ R_n(Q)$ is dense in $X^n$.
\item There exists a dense subset $A$ of $ X$ which is a
$\{ R_n(Q) \}$ set, i.e. $A^n \subset R_n(Q)$ for $n = 1,2, \dots $.
\item $Q$ is dense in $C'(X)$.
\item $CANTOR(X) \cap Q$ is a dense $G_{\delta}$ subset of $C'(X)$.
\item There is a sequence $\{ K_i : i = 1,2, \dots \}$ which is dense in $CANTOR(X)$ such that $\bigcup_{i=1}^n \ K_i  \ \in \ Q$ for
$n = 1,2, \dots$.
\end{enumerate}

(b) The following conditions are equivalent
\begin{enumerate}
\item There is a Cantor set in $Q$, i.e. $CANTOR(X) \cap Q \not = \emptyset$.
\item There is a Cantor set which is an $\{ R_n(Q) \}$ set.
\item There is an uncountable  $\{ R_n(Q) \}$ set.
\item There is a nonempty  $\{ R_n(Q) \}$ set with no isolated points.
\item There is a nonempty, closed, perfect subset $Y$ of $X$ such that
$Y^n \cap R_n(Q)$ is dense in $Y^n$ for $n = 1,2, \dots$.
\end{enumerate}
\end{thm}

\medskip

\subsection{Uniform chaos in light of the Kuratowski-Mycielski Theorem}

With this new vocabulary we can restate
Theorem \ref{tran} by saying that for a transitive system $(X,T)$ the collection $Q(Recur)$, of uniformly rigid subset, is a dense $G_\delta$ subset of $C'(X)$.
For the reader's convenience we repeat the statement of the theorem
(augmented with a statement about pairwise proximality) and provide
a short proof which employs the Kuratowski-Mycielski machinery.

\begin{thm} \label{tran'}
Let $(X,T)$ be a transitive TDS without isolated points. Then
there are Cantor sets $C_1 \subseteq C_2 \subseteq \cdots$ such
that $\bigcup \limits_{i=1}^{\infty}C_n$ is a dense rigid subset
of $X_{tr}$ and for each $N\in \N$, $C_N$ is uniformly rigid.
\begin{itemize}
\item
If in addition, $P(X,T)$ is dense in $X \times X$ then we can
require that $\bigcup \limits_{i=1}^{\infty}C_n$ is pairwise proximal.
\item
If in addition, for each $n\in \N$, $Prox_n(X)$ is dense in
$X^n$, then we can require that for each $N\in \N$, $ C_N$ is
uniformly  proximal. Thus under these conditions $(X,T)$ is uniformly chaotic.
\end{itemize}
\end{thm}

\begin{proof}
As described in Example (1) above, $Recur_n$ and $(X_{tr})^n$ are dense $G_{\delta}$ subsets of $X^n$.
Hence, condition (1) of  part (a) of the Kuratowski-Mycielski Theorem
applies to $Q(Recur) \cap Q(X_{tr})$.
The result follows from condition (5) of part (a) with
$C_N = \bigcup_{i=1}^N \ K_i$.

If $P(X,T)$ is dense in $X^2$ then we can intersect as well with the dense $G_{\delta}$ set $Q(P(X,T))$.

If $Prox_n$ is dense in $X^n$ for every $n$ then  $Q(Prox)$ is also a dense $G_{\delta}$ by the Kuratowski-Mycielski Theorem and so we can intersect
with it as well.
\end{proof}

\begin{rem}
Notice that in general the collection $Q(Recur)$ of uniformly rigid
subsets of $X$, is not {\em finitely determined};
that is, a closed subset
$A \subset X$ with $A^n \subset Recur_n$ for every $n \ge 1$ is
merely weakly rigid and need not be uniformly rigid.
Similarly $Q(Prox)$ is not finitely determined and
a closed subset $A \subset X$ with $A^n \subset Prox_n$ for every $n \ge 1$ is
merely a proximal set and need not be uniformly  proximal.
\end{rem}

\medskip

%
%

We do likewise with the criterion for chaos (Theorem \ref{zhunzhe}).

\begin{thm}[A criterion for chaos] \label{zhunzhe'}
Let $(X,T)$ be a transitive TDS without isolated points. Assume that
$(Y,T)$ is a subsystem of $(X,T)$ such that $(X\times Y,T)$ is
transitive,
there are Cantor sets $C_1 \subseteq C_2 \subseteq
\cdots$ such that
\begin{enumerate}
\item
$K = \bigcup \limits_{i=1}^{\infty}C_n$ is a dense subset of
$X_{tr}$ and;
\item
for each $N\in \N$, $C_N$ is
a Kronecker set for $Y$ and is uniformly rigid.
\end{enumerate}

In particular, $(X,T)$ is densely uniformly  chaotic.
\end{thm}

\begin{proof}
We follow the notation of Examples (1) and (3) above.
The work below will be to show that $R_n(Q(TRANS,Y))$ is dense in $X^n$ for $n = 1,2, \dots$.  We have already seen that $Recur_n$ is dense in $X^n$.
By the Kuratowski-Mycielski Theorem  it follows that
$$
Q(TRANS,Y) \cap Q(Recur) \cap Q(X_{tr})
$$
is dense in $C'(X)$ and that the required sequence of Cantor sets exists.

Fix $\ep > 0$ and $y_1, \dots ,y_n \in Y$ and choose  open subsets
$W_1, \dots ,W_n$  of diameter less than $\ep$
with $y_i \in W_i$ for $i = 1, \dots ,n$.
We will prove that the open set
$\bigcup_{k \in \Z_+} \ T^{-k}W_1 \times \cdots \times T^{-k}W_n$ is  dense.
Then intersect over positive rational $\ep$ and
$\{y_1, \dots ,y_n\}$ chosen from a countable dense subset of $Y$.
The Baire Category Theorem then implies that
$R_n(Q(TRANS,Y))$ is a dense $G_{\delta}$ subset of $X^n$ as required.

Let $U_1, \dots ,U_n$ be open nonempty subsets of $X$.
Because $X \times Y$ is  transitive there exists
$r_2 \in N(U_1 \times (W_1 \cap Y), U_2 \times (W_2 \cap Y))$.
Let
$$
U_{12} \times W_{12} =
(U_1 \cap T^{-r_2}U_2)\times (W_1 \cap T^{-r_2}W_2),
$$
an open set which meets $X \times Y$.
Proceed inductively, finally choosing
$r_n \in N((U_{1 \dots n-1} \times (W_{1 \dots n-1}\cap Y)),
U_n \times (W_n \cap Y))$
and let
$$
U_{1 \dots n} \times W_{1 \dots n}
= (U_{1 \dots n-1} \cap T^{-r_n} U_n)
\times (W_{1 \dots n-1} \cap T^{-r_n}W_n).
$$

Choose $(x,y) \in (U_{1\dots n} \times W_{1 \dots n})
\cap (Z \times Y)$ with $x \in X_{tr}$.
Thus, $(x,T^{r_2}x, \dots ,T^{r_n}x) \in U_1 \times \cdots \times U_n$ and
$(y,T^{r_2}y, \dots ,T^{r_n}y) \in W_1 \times \cdots \times W_n$.
Since $x$ is a transitive point, we
can choose $T^kx$ close enough to $y$ so that
$(T^kx,T^{k+r_2}x, \dots ,T^{k+r_n}x) \in W_1 \times \cdots \times W_n$.
Thus, $(x,T^{r_2}x, \dots ,T^{r_n}x) \in (U_1 \times \cdots \times U_n)
\cap  (T^{-k}W_1 \times \cdots \times T^{-k}W_n)$, as required.

For the last assertion of the theorem use Lemma \ref{Kronecker}.
\end{proof}

\medskip

\section{Chaotic subsets of minimal systems}\label{Sec-min}
It is well known that a non-equicontinuous minimal system is sensitive
(see \cite{GW93}).
In this section we will have a closer look at
chaotic behavior of minimal systems and will examine the relationship
between chaos and structure theory.

\subsection{On the structure of minimal systems}
The structure theory of minimal systems originated
in Furstenberg's seminal work \cite{F63}.
In this subsection we briefly review some of the main
results of this theory. It was mainly developed for
group actions and accordingly we assume for the rest of
the paper that $T$ is a homeomorphism.
Much of this work can be done for a general locally
compact group actions, but for simplicity we stick to the traditional
case of $\Z$-actions. We refer the reader to \cite{G76}, \cite{V77}, 
and \cite{Au} for details.

We first recall that an extension
$\pi : X \to Y$ of minimal systems is called a {\em relatively
incontractible (RIC) extension}\ if it is open and for every $n \ge 1$
the minimal points are dense in the relation
$$
R^n_\pi = \{(x_1,\dots,x_n) \in X^n : \pi(x_i)=\pi(x_j),\ \forall \ 1\le i
\le j \le n\}.
$$
(See Theorem \ref{RIC} in the appendix below.)

We say that a minimal system $(X,T)$ is a
{\em strictly PI system} if there is an ordinal $\eta$
(which is countable when $X$ is metrizable)
and a family of systems
$\{(W_\iota,w_\iota)\}_{\iota\le\eta}$
such that (i) $W_0$ is the trivial system,
(ii) for every $\iota<\eta$ there exists a homomorphism
$\phi_\iota:W_{\iota+1}\to W_\iota$ which is
either proximal or equicontinuous
(isometric when $X$ is metrizable), (iii) for a
limit ordinal $\nu\le\eta$ the system $W_\nu$
is the inverse limit of the systems
$\{W_\iota\}_{\iota<\nu}$,  and
(iv) $W_\eta=X$.
We say that $(X,T)$ is a {\em PI-system} if there
exists a strictly PI system $\tilde X$ and a
proximal homomorphism $\theta:\tilde X\to X$.

If in the definition of PI-systems we replace
proximal extensions by almost one-to-one
extensions (or by highly proximal extensions
in the non-metric case) we get the notion of HPI {\em systems}.
If we replace the proximal extensions by trivial
extensions (i.e.\ we do not allow proximal
extensions at all) we have I {\em systems}.
These notions can be easily relativize and we then speak
about I, HPI, and PI extensions.

In this terminology Furstenberg's structure theorem for distal
systems (Furstenberg \cite{F63}) and
the Veech-Ellis structure theorem for point distal systems
(Veech \cite{V70}, and Ellis \cite{E}),
can be stated as follows:

\begin{thm}\label{FST}
A metric minimal system is distal if and only if  it is an I-system.
\end{thm}

\begin{thm}
A metric minimal dynamical system is point distal if and only if  it is an
HPI-system.
\end{thm}

Finally we have the structure theorem for minimal systems,
which we will state in its relative form
(Ellis-Glasner-Shapiro \cite{EGS},
McMahon \cite{Mc}, Veech \cite{V77}, and
Glasner \cite{G}).

\begin{thm}[Structure theorem for minimal systems]\label{structure}
Given a homomorphism $\pi: X \to Y$ of minimal dynamical system,
there exists an ordinal $\eta$
(countable when $X$ is metrizable) and a canonically defined
commutative diagram (the canonical PI-Tower)
\begin{equation*}
\xymatrix
        {X \ar[d]_{\pi}             &
     X_0 \ar[l]_{{\theta}^*_0}
         \ar[d]_{\pi_0}
         \ar[dr]^{\sigma_1}         & &
     X_1 \ar[ll]_{{\theta}^*_1}
         \ar[d]_{\pi_1}
         \ar@{}[r]|{\cdots}         &
     X_{\nu}
         \ar[d]_{\pi_{\nu}}
         \ar[dr]^{\sigma_{\nu+1}}       & &
     X_{\nu+1}
         \ar[d]_{\pi_{\nu+1}}
         \ar[ll]_{{\theta}^*_{\nu+1}}
         \ar@{}[r]|{\cdots}         &
     X_{\eta}=X_{\infty}
         \ar[d]_{\pi_{\infty}}          \\
        Y                 &
     Y_0 \ar[l]^{\theta_0}          &
     Z_1 \ar[l]^{\rho_1}            &
     Y_1 \ar[l]^{\theta_1}
         \ar@{}[r]|{\cdots}         &
     Y_{\nu}                &
     Z_{\nu+1}
         \ar[l]^{\rho_{\nu+1}}          &
     Y_{\nu+1}
         \ar[l]^{\theta_{\nu+1}}
         \ar@{}[r]|{\cdots}         &
     Y_{\eta}=Y_{\infty}
    }
\end{equation*}
where for each $\nu\le\eta, \pi_{\nu}$
is RIC, $\rho_{\nu}$ is isometric, $\theta_{\nu},
{\theta}^*_{\nu}$ are proximal and
$\pi_{\infty}$ is RIC and weakly mixing of all orders.
For a limit ordinal
$\nu ,\  X_{\nu}, Y_{\nu}, \pi_{\nu}$
etc. are the inverse limits (or joins) of
$ X_{\iota}, Y_{\iota}, \pi_{\iota}$ etc. for $\iota
< \nu$.
Thus $X_\infty$ is a proximal extension of $X$ and a RIC
weakly mixing extension of the strictly PI-system $Y_\infty$.
The homomorphism $\pi_\infty$ is an isomorphism (so that
$X_\infty=Y_\infty$) if and only if  $X$ is a PI-system.
\end{thm}

\medskip

\subsection{Lifting chaotic sets}

Using enveloping semigroup techniques we are able
to lift chaotic sets in minimal systems. First a lemma
concerning proximal sets in minimal systems.


\begin{lem}\label{cell}
Let $(X,T)$ be a minimal TDS and $K \subset X$ a proximal set. The set
$$
A_K=\{p\in \b^* \Z: pK\ {\text{is a singleton}}\}
$$
is a closed ideal in $\b^* \Z$.  If $I \subset A_K$ is any minimal ideal
then for any $x_0 \in K$,
$$
K \subset \{vx_0: v \in J(I)\}.
$$
\end{lem}

\begin{proof}
Because
$K$ is a proximal set, $A_K$ is nonempty and is clearly a left ideal. Since
$$
A_K=\{p\in \b^* \Z : px_1=px_2 \ {\text{for all}}\  x_1, x_2 \in K \}
$$
it is closed as well. Let $I \subset A_K$ be a minimal left ideal,
which exists by Ellis' theory.  Because $(X,T)$ is minimal, $Ix = X$ for
any $x \in K$ and so $I_x = \{p\in I: px = x \}$ is a nonempty closed
subsemigroup. By Ellis' Lemma there exists an idempotent $v_x \in I_x$ and
since $v_x \in A_K$ and $x \in K$ we have $v_xK = \{ x \}$.
Thus, for any $x_0 \in K$
$$
K = \{v_xx_0 : x \in K \}.
$$
\end{proof}

\begin{rem}
Because $I$ is a minimal ideal, $I = Iu$ for any idempotent 
$u \in I$ and so $pu = p$ for any $p \in I$.  In particular,
$$
v_{y}v_{x} = v_{y} \qquad \text{for all} \ x,y \in K.
$$
\end{rem}

%

%

\medskip

\begin{lem}\label{lift}
Let $\pi:X\lra Y$ be an extension between minimal systems.
\begin{enumerate}
\item
If $\pi$ is a proximal extension and
$K\subset Y$ is a  proximal set of $Y$,
then any set $K'$ of $X$ with $\pi(K')=K$ is a proximal set.
\item
For any proximal subset $K$ of $Y$ there is
a proximal subset $K'$ of $X$ with $\pi(K')=K$.
\item
For any weakly rigid subset $K$ of $Y$, there is
a weakly rigid subset $K'$ of $X$ with $\pi(K')=K$.
Moreover if $K$ is both proximal and weakly rigid then
there is a subset $K'$ of $X$ with $\pi(K')=K$
which is both proximal and weakly rigid.
In particular, for any strongly Li-Yorke pair $(y,y')$ in $Y\times Y$ there is
a strongly Li-Yorke pair $(x,x')$ in $X\times X$
with $\pi(x)=y, \pi(x')=y'$.
\item
If $\pi$ is a distal extension and $K\subset Y$ is a weakly rigid
set of $Y$, then any set $K'$ of $X$ with $\pi(K')=K$ is a weakly
rigid set.
\end{enumerate}
In the cases {\rm (2)} and {\rm (3)} we have $\pi \rest K'$ is one-to-one.
\end{lem}

\begin{proof}  If $K$ is a proximal subset of $Y$ we apply Lemma 
\ref{cell} and its proof to define the
ideal $A_K$ in $\b^* \Z$, choose a minimal ideal $I \subset A_K$ 
and idempotents $\{ v_x \in I : x \in K \}$
such that $v_xx = x$ for all $x \in K$.

1.\
%
%
Now assume that $\pi$ is proximal and $\pi(K') = K$ with $K$ a proximal subset.
Let $u $ be an arbitrary idempotent in $I$ so that $uK$ is a singleton.

For any pair
$x_1,x_2\in K'$  we have
$$\pi(ux_1)=u\pi(x_1)=u\pi(x_2)=\pi(ux_2).$$

As $u(ux_1,ux_2)= (ux_1,ux_2)$ and $u$ is a minimal idempotent,
$(ux_1,ux_2)$ is a minimal point. Since $\pi$ is proximal, we have
$ux_1=ux_2$. Since the pair $x_1,x_2$ was arbitrary, $uK'$ is a singleton.

\medskip

2.\
Fix $x_0 \in X$ such that $y_0 = \pi(x_0) \in K$. 
Assuming that $K$ is a proximal subset we define
$j : K \to X$ by $j(x) = v_xx_0$. 
Observe that $\pi(j(x)) = \pi(v_xx_0) = v_xy_0 = x$. So with
$K' = j(K)$ we have $\pi(K') = K$. On the other hand, 
$y \in K$ implies $v_yv_x = v_y$ and so $v_yj(x) = v_yx_0 = j(y)$ 
for all $x \in K$. That is, $v_yK'$ is the singleton $\{ j(y) \}$ and so
$K'$ is proximal.

%
%
%

\medskip

3.\
Assume that $K$ is a weakly rigid subset. The set
$$
S_K = \{p \in \b^* \Z : py = y\ \text{for every}\ y \in K \}
$$
is a closed subsemigroup, nonempty because $K$ is weakly rigid. 
By Ellis' Lemma there is an idempotent $u \in S_K$. 
Choose for each $x \in K$, $j(x) \in \pi^{-1}(x)$.
Let 
$$ 
K' = \{ uj(x) : x \in K \}.
$$ 
Since $\pi(uj(x)) = u\pi(j(x)) = ux = x$ it follows that
$\pi(K') = K$. Since $u$ is an idempotent it acts as the identity on $K'$.

Now assume in addition that $K$ is a proximal subset. 
$A_Ku$ is a closed ideal. Since $u$ acts as the
identity on $K$, it follows that $pu(K)$ is a singleton for every 
$p \in A_K$, i.e. $A_Ku \subset A_K$.
If $I$ a minimal ideal in $A_Ku$ then with $pu = p$ for all $p \in I$.  
In particular, the idempotents $v_x \in I$ satisfy $v_x u = v_x$ and so 
$u v_xuv_x = u v_xv_x = uv_x$. That is, $uv_x $ is an idempotent in
$I$. Furthermore, $uv_x(K) = \{ ux \} = \{ x \}$. Thus, we can replace $v_x$ by $uv_x$ if necessary and so assume that $uv_x = v_x$.

As in (2) define $j(x) = v_xx_0$ to obtain the proximal set $K' = j(K)$. Since $uv_x = v_x$, $uj(x) = j(x)$ and
so $u$ acts as the identity on $K'$.  That is, $K'$ is a weakly rigid set as well.




\medskip

4.\ 
As in part 3. set
$$
S_K = \{p \in \b^* \Z : py = y\ \text{for every}\ y \in K \}
$$
and then pick an idempotent $u \in S_K$.
Now for any $x \in X$ with $\pi(x) \in K$ the points
$x$ and $ux$ are proximal. But as $\pi(ux) =
u\pi(x) =\pi(x)$ and $\pi$ is a distal extension we conclude
that $ux=x$. Thus if $\pi(K')=K$ then $ux=x$ for every
$x \in K'$, whence $K'$ is weakly rigid. 
\end{proof}

Now from Theorem \ref{KMT}
it follows that a dynamical system $(X,T)$ contains
a Cantor subset which is both uniformly  proximal
and uniformly rigid, if and only if there is an uncountable
$A\subset X$ such that for every $n$-tuple $(x_1,\dots,x_n)$ with $x_j \in A$
we have $(x_1,\dots,x_n) \in Prox_n(X) \cap Recur_n(X)$.
Thus if we let $Q(PR)=Q(Prox)\cap Q(Recur)$
be the collection of closed subsets of
$X$ which are both uniformly  proximal and uniformly rigid,
then $Q(PR)$ is a $G_\delta$ subset of $C'(X)$ and for every $n\ge 1$,
$R_n(Q(PR)) = Prox_n(X) \cap Recur_n(X)$.
These facts combined with Lemma \ref{lift}(3) yield the following important corollary.

\begin{thm}\label{Ak}
Let $\pi: X \to Y$ be a homomorphism of minimal systems.
If $Y$ contains a uniformly  chaotic subset then so does $X$.
\end{thm}

\begin{rem}
One would like to prove analogous lifting theorems for
Li-Yorke and strong Li-Yorke chaotic sets
(i.e. uncountable scrambled and strongly scrambled sets).
Unfortunately the collection of closed scrambled sets
is not, in general, a $G_\delta$ subset of $C'(X)$,
and we therefore can not use this kind of argument to show
that Li-Yorke chaos lifts under homomorphisms of minimal
systems. The problem with lifting closed strongly scrambled sets
(which do form a $G_\delta$ set) is that we do not know
whether an uncountable strongly scrambled set can always
be lifted through an extension of minimal systems.
\end{rem}

\medskip

\subsection{Weakly mixing extensions}






\begin{thm} \label{wm}
Let $(X,T)$ be a TDS and $\pi: (X,T)\rightarrow (Y,S)$ an open
nontrivial weakly mixing extension. Then there is a residual
subset $Y_0\subseteq Y$ such that for every point $y\in Y_0$ the
set $\pi^{-1}(y)$ contains a dense strongly scrambled Mycielski subset
$K$ such that $K\times K \setminus \Delta_X \subseteq Trans(R_{\pi})$.
In particular $(X,T)$ is strongly Li-Yorke chaotic.

If moreover $\pi$ is weakly mixing and RIC, then there is a residual
subset $Y_0\subseteq Y$ such that for every point $y\in Y_0$ a
dense Mycielski set $K \subset \pi^{-1}(y),\ y\in Y_0$ as above can be found which is uniformly chaotic, whence $X$ is uniformly chaotic.
%

\end{thm}

\begin{proof}
Since $\pi$ is open, it follows that $\pi \times \pi: R_{\pi}\to Y,
(x_1,x_2)\mapsto \pi(x_1)$ is open as well. Since $R_{\pi}$ is
transitive, the set of transitive points $Trans(R_{\pi})$ is a
dense $G_{\delta}$ subset of $R_{\pi}$. By Ulam's Theorem there is a
residual subset $Y_0\subseteq Y$ such that for every point $y\in
Y_0$,
$$
Trans(R_{\pi}) \cap Recur_2 \cap (\pi_\infty^{-1}(y)\times
\pi_\infty^{-1}(y))
$$
is dense $G_{\delta}$ in $(\pi^{-1}(y)\times \pi^{-1}(y))$.

Now for each $y\in Y_0$, we claim that $\pi^{-1}(y)$ has no isolated
points. In fact if this is not true, then there exists $x\in
\pi^{-1}(y)$ such that $\{ x \}$ is an open subset of $\pi^{-1}(y)$.
Moreover, $\{ (x,x) \}$ is an open subset of $\pi^{-1}(y)\times
\pi^{-1}(y)$. Since $Trans(R_{\pi}) \cap (\pi^{-1}(y)\times
\pi^{-1}(y))$ is dense $G_{\delta}$ in $(\pi^{-1}(y)\times
\pi^{-1}(y))$, one has $(x,x)\in Trans(R_{\pi})$. This shows that
$R_{\pi}=\Delta_X$ which contradicts the fact that $\pi$ is a
non-trivial extension. Finally, by Theorem \ref{special}
there is a dense
s-chaotic subset $K\subseteq \pi^{-1}(y)$ such that
\begin{equation*}
K\times K \setminus \Delta_X \subseteq Trans(R_{\pi}) \cap
(\pi^{-1}(y)\times \pi^{-1}(y))\setminus \Delta_X \subseteq
Trans(R_{\pi}).
\end{equation*}

We now further assume that $\pi$ is a RIC extension. Then
by \cite{G}, Theorem 2.7, $\pi$ is weakly mixing of all orders
(i.e. $R^n_\pi$ is transitive for all $n \ge 2$) and in particular
for every $n \ge 2$,
$$
Prox_n \cap Recur_n \cap \pi^{-1}(y)^n
$$
is a dense in $\pi^{-1}(y)^n$, for every $y \in Y_0$.
Applying Theorem \ref{KMT} we obtain our claim.

%
\end{proof}

\medskip

\subsection{The non $PI$ case}\label{non-PI}

%

%



The following theorems of Bronstein \cite{Bro} and van der Woude \cite{Wo85}
give intrinsic characterizations of PI-extensions and HPI-extensions
respectively.
Recall that a map $\pi : X \to Y$ between compact spaces is called
{\em semi-open} if ${\rm int}\,\pi(U) \ne \emptyset$ for every nonempty
open subset $U \subset X$. It was observed by J. Auslander
and N. Markley that a homomorphism $\pi : X \to Y$
between minimal systems is always semi-open (see e.g. \cite{G},
Lemma 5.3).

\begin{thm}\label{bw}
Let $\pi:(X,T)\lra (Y,T)$ be a homomorphism of compact metric
minimal systems. Then
\begin{enumerate}
\item
The extension $\pi$ is PI if and only if it satisfies the
following property: whenever $W$ is a closed invariant subset of
$R_\pi$ which is transitive and has a dense subset of
minimal points, then $W$ is minimal.
\item
The extension $\pi$ is HPI if and only if it satisfies the
following property: whenever $W$ is a closed invariant subset of
$R_\pi$ which is transitive and the restriction of
the projection maps to $W$ are semi-open, then $W$ is minimal.
\end{enumerate}
\end{thm}

Next we show that a minimal system which is
a non-PI extension has an s-chaotic subset.

\begin{thm}\label{nonpi}
Let $\pi:(X,T)\lra (Y,T)$ be a homomorphism of
metric minimal systems. If $\pi$ is a non-PI extension, then there is a
dense subset $Y_0\subset Y$ such that for each $y_0\in Y_0$, there
is a uniformly chaotic subset of $\pi^{-1}(y_0)$.
In particular $(X,T)$ is strongly Li-Yorke chaotic.
\end{thm}

\begin{proof}
Assume that $\pi:(X,T)\lra (Y,T)$ is a non-PI extension.
Then by Theorem \ref{structure} there exist
$\phi:(X_{\infty},T)\rightarrow (X,T)$,
$\pi_{\infty}:(X_{\infty},T)\rightarrow (Y_{\infty},T)$ and
$\eta:Y_\infty\lra Y$ such that $\phi$ is a proximal extension,
$\pi_{\infty}$ is weakly mixing RIC extension, and $\eta$
is a PI-extension. As $\pi$ is non-PI, $\pi_\infty$ is
non-trivial.

Now consider the commutative diagram
\begin{equation*}
\xymatrix
{
X \ar[d]_{\pi}  &  X_\infty \ar[l]_{\phi}\ar[d]^{\pi_\infty} \\
Y &  Y_\infty\ar[l]^{\eta}
}
\end{equation*}


By Theorem \ref{wm} there is a dense $G_\delta$ subset $Y^0_\infty
\subset Y_\infty$ such that, for every $y \in Y^0_\infty$,
there is a dense uniformly chaotic subset $K_y$ of $\pi^{-1}_\infty(y)$.

Since $\pi$ is not PI, $\pi_\infty$ is not proximal. Thus, there
is a distal point $(x_1,x_2)\in R_{\pi_\infty}\setminus
\Delta_{X_\infty}$. This implies that $\phi(x_1)\neq \phi(x_2)$ as
$\phi$ is a proximal extension.
For any $k_1,k_2\in K = K_y$ with $k_1\neq k_2$, one has $(k_1,k_2)\in
Trans(R_{\pi_\infty})$. As $\phi(x_1)\neq \phi(x_2)$, $(x_1,x_2)
\in R_{\pi_\infty}$ and $(k_1,k_2)\in Trans(R_{\pi_\infty})$, one
has $\phi(k_1)\neq \phi(k_2)$. That is, $\phi:K\rightarrow
\phi(K)$ is a bijection. Therefore, as is easy to check, $\phi(K)$ is
a uniformly chaotic subset of $X$. Moreover $\phi(K)\subset
\pi^{-1}(\eta(y))$. Finally we let $Y_0=\eta(Y^0_\infty)$; clearly
a dense subset of $Y$.
\end{proof}

\medskip

\subsection{The proximal but not almost one-to-one case}

Every extension of minimal systems can be lifted to an open
extension by almost one-to-one modifications. To be precise, for
every extension $\pi:X\rightarrow Y$ of minimal systems there
exists a canonically defined commutative diagram of extensions
(called the {\em shadow diagram})
\begin{equation*}
\xymatrix
{
X \ar[d]_{\pi}  &  X^* \ar[l]_{\sigma}\ar[d]^{\pi*} \\
Y &  Y^*\ar[l]^{\tau}
}
\end{equation*}

with the following properties:
\begin{enumerate}
\item[(a)]
$\sigma$ and $\tau$ are almost one-to-one;
\item[(b)]
$\pi^*$ is an open extension;
\item[(c)]
$X^*$ is the unique minimal set in $R_{\pi
\tau}=\{(x,y)\in X\times Y^*: \pi(x)=\tau (y)\}$ and $\sigma$ and
$\pi^*$ are the restrictions to $X^*$ of the projections of
$X\times Y^*$ onto $X$ and $Y^*$ respectively.
\end{enumerate}

\medskip

In \cite{G76} it was shown that a metric minimal system $(X,T)$ with the 
property that $Prox_n(X)$ is dense in $X^n$ for every $n \ge 2$ is weakly 
mixing. This was extended by van der Woude  \cite{Wo82} as follows (see 
also \cite{G}).

\begin{thm}\label{pwm}
Let $\pi: X \to Y$ be a factor map of the metric minimal system $(X,T)$.
Suppose that $\pi$ is open and that for every $n \ge 2$,
$Prox_n(X)\cap R_\pi$ is dense in $R_\pi$. Then $\pi$ is a weakly 
mixing extension. In particular a nontrivial open proximal extension is
a weakly mixing extension.
\end{thm}

\begin{lem}\label{1-1}
Let $\pi:X\rightarrow Y$ be a continuous surjective map between compact
metric spaces which is almost one-to-one. If $A\subset X$ is a dense
$G_\delta$ subset, then $\pi(A)$ contains a dense $G_\delta$ subset of $Y$.
\end{lem}
\begin{proof}
Let $A_0=\{ x\in X: \pi^{-1}\pi(x)=\{ x \}\}$ and $B_0=\{ y\in
Y:\text{Card}\ \pi^{-1}(y)=1 \}$. Then $A_0$ (resp. $B_0$) is a
dense $G_\delta$ subset of $X$ (resp. $Y$). Now $A\cap A_0$ is a
dense $G_\delta$ subset of $X$, hence a dense
$G_\delta$ of $A_0$. As the set of continuity points of
$\pi^{-1}:Y\lra C(X)$ contains $B_0$, $\pi:A_0\rightarrow B_0$ is a
homeomorphism, so $\pi(A\cap A_0)$ is a dense $G_\delta$ subset of
$B_0$. Therefore, there exist open subsets $U_n$ of $Y$ such that
$\bigcap_{n=1}^\infty U_n \cap B_0=\pi(A\cap A_0)$. This shows
that $\pi(A\cap A_0)$ is also a dense $G_\delta$ subset of $Y$.
\end{proof}

\medskip

Recall that a subset $K$ of $X$ is a proximal set
if each finite tuple from $K$ is uniformly  proximal (see Definition
\ref{proximal}). The proof of the following lemma is straightforward.

\begin{lem}
Let $\pi:X\lra Y$ be a proximal extension between minimal systems. Then
for each $y\in Y$, $\pi^{-1}(y)$ is a proximal set.
\end{lem}

In the sequel it will be convenient to have the following:

\begin{de}\label{chaotic}
Let $(X,T)$ be a TDS.
\begin{enumerate}
\item
A scrambled Mycielski subset $K \subset X$ will be called a
{\em chaotic subset} of $X$.
\item
A strongly scrambled Mycielski subset will be called an
{\em s-chaotic subset} of $X$.
\end{enumerate}
\end{de}

We can now prove the following result (see also \cite{AAG}, Theorem 6.33).

\begin{thm}\label{pn1-1}
Let $\pi:X\rightarrow Y$ be a proximal but not almost one-to-one
extension between minimal systems.
Then there is a residual subset $Y_0\subset Y$ such
that for each $y\in Y_0$,  $\pi^{-1}(y)$ contains a proximal s-chaotic set
$K$ with $K\times K \setminus \Delta_X \subseteq Trans(R_{\pi})$.
\end{thm}
\begin{proof}
In the shadow diagram for $\pi$, the map $\pi^*$ is
open and proximal. Since $\pi$ is not almost one-to-one $\pi^*$ is
not trivial. 
Thus, by Theorem \ref{pwm}, $\pi^*$ is a nontrivial open weakly 
mixing extension.
Hence by Theorem \ref{wm} there is
a residual subset $Y^*_0\subset Y^*$ such that for each $y^*\in
Y^*_0$, $\pi^{*-1}(y^*)$ contains an s-chaotic set $K^*$ and
$K^*\times K^* \setminus \Delta_{X^*} \subseteq
Trans(R_{\pi^*})$. Moreover, $\pi^*$ being proximal, we have for
every $n \ge 2$,
$\pi^{*-1}(y^*)^n \subset Prox_n$ and therefore we can require that $K^*$
be proximal as well.
Since in the shadow diagram $\sigma$ and
$\pi^*$ are the restrictions to $X^*$ of the projections of
$X\times Y^*$ onto $X$ and $Y^*$ respectively, $\sigma (K^*)$ is an
s-chaotic set, as $\pi\sigma(K^*)
=\tau\pi^*(K^*)=\{\tau(y^*)\}$, $\sigma(K^*)\subset
\pi^{-1}(\tau(y^*))$.
Finally, set $Y_0=\tau(Y^*_0)$. Since $\tau$ is almost one-to-one, $Y_0$ is a
residual subset of $Y$ (Lemma \ref{1-1}) .
\end{proof}

The following result was first proved in \cite{AAG}.

\begin{cor}
Let $\pi:X\rightarrow Y$ be an asymptotic extension between
minimal systems. Then $\pi$ is almost one-to-one.
\end{cor}
\begin{proof} We use the notations in the proof of Theorem \ref{pn1-1}.
Note that $\pi$ is proximal. If it is not almost one-to-one, then
by Theorem \ref{pn1-1}, there are $\sigma(x)\not=\sigma(y)\in
\sigma(K^*)$ such that $(\sigma(x),\sigma(y))$ is a recurrent
point of $T\times T$. It is clear that $\pi\sigma(x)=\pi\sigma(y)$
since the diagram is commutative. It now follows that $\pi$ is not
asymptotic, a contradiction.
\end{proof}

\begin{rem}\label{w}
In \cite{GW79} the authors construct an example of a
minimal system $(X,T)$ which admits a factor map
$\pi: X \to Y$ such that (i) the factor $Y$ is equicontinuous,
(ii) the map $\pi$ is a nontrivial open proximal extension.
Now such an $X$ is clearly strictly PI but not HPI.
However, according to Theorem \ref{pwm} the extension $\pi$
is a weakly mixing extension and it follows from Theorem \ref{wm}
that for some $y \in Y$ the fiber $\pi^{-1}(y)$ contains a
dense proximal s-chaotic subset. 
Thus $X$ is an example of a minimal PI system which is strongly Li-Yorke chaotic. We do not have an example of a minimal PI system which 
contains a uniformly chaotic set.  

We also note that, by \cite{BGKM}, positive topological entropy
implies the existence of an s-chaotic subset. Since there are 
HPI systems with positive entropy (e.g many Toeplitz systems
\cite{W}) we conclude that there are HPI systems which are strongly
Li-Yorke chaotic.

\end{rem}

\medskip

\subsection{The $PI$, non-$HPI$ case}
We have shown (Subsection \ref{non-PI}) that for a non-PI extension
there is a uniformly chaotic set.
The natural question now is whether there is a chaotic
(s-chaotic, uniformly chaotic) set for a non-HPI extension?
(Recall that for a metric $X$ the notions  `HPI extension' and
`point distal extension' coincide.)
At present we are unable to answer this question fully.
However we will show that the answer is affirmative
for a sub-class of non-HPI extensions:

\begin{prop}\label{spi}
Let $\pi:X\lra Y$ be a strictly PI extension but non-HPI extension
between minimal systems. Then there is a dense set $Y_0$ of $Y$
such that for each $y\in Y_0$, $\pi^{-1}(y)$ contains a proximal
chaotic set.
\end{prop}
\begin{proof} Since by assumption $\pi$ is non-HPI, in its strictly PI-tower
at least one of the proximal extensions in the canonical PI-tower
is not an almost one-to-one extension. Let us denote this segment of
the tower by
$$
X \overset{\pi_1}\longrightarrow Z_1\overset{\pi_2}\longrightarrow Z_2
\overset{\pi_3}\longrightarrow Y,
$$
with $\pi_1\circ \pi_2\circ \pi_3=\pi$, and where
$\pi_1$ and $\pi_3$ are strictly PI extensions and
$\pi_2$ is a proximal but not an almost one-to-one extension.

By Theorem \ref{pn1-1} there exists a dense set $Z_0\subset Z_2$
such that for each $z\in Z_2$,
$\pi_2^{-1}(z)$ contains a proximal s-chaotic set
$K \subset Z_1$ .  By  Lemma \ref{lift} (3), there is a proximal
subset $K' \subset X$ with $\pi_1 \circ \pi_1(K') =K$
and as $K$ is s-scarambled, $K'$ is at least scrambled.
Now the proposition follows by setting $Y_0=\pi_3(Z_0)$.
\end{proof}

Now assume that $\pi:X \lra Y$ is PI and not HPI. This means that
in the canonical PI-tower there are maps $\phi:X_\infty\lra X$ which
is proximal and $\eta: X_\infty\lra Y$ which is strictly PI:
\begin{equation*}
\xymatrix
{
& X_\infty\ar[dl]_{\phi} \ar[dr]^{\eta}  &  \\
X \ar[rr]^\pi & & Y
}
\end{equation*}

\medskip

\begin{lem}\label{case1}
The extension $\eta$ is not strictly HPI
extension.
\end{lem}
\begin{proof}
Assume $\eta$ is a strictly HPI-extension. Then by Theorem \ref{bw}.(2)
there are no non-minimal transitive subsystem $W$ of
$R_{\eta}$ such that the coordinate projection $W\rightarrow Y$ is
semi-open. Now it is easy to see that there is no non-minimal
transitive subsystem $W'$ of $R_{\phi}$ such that the coordinate
projection $W'\rightarrow X$ is semi-open. For if $W'$ is a
non-minimal transitive subsystem of $R_{\phi}$ such that the
coordinate projection $W'\rightarrow X$ is semi-open, then $W'$ is
also a subsystem of $R_{\eta}$. But the composition of two
semi-open maps is also semi-open and $\pi:X\rightarrow Y$ is
semi-open, hence the coordinate projection $W'\rightarrow Y$ is
semi-open, a contradiction.

Hence using Theorem \ref{bw}.(2) again, this shows that $\phi$
is an HPI-extension. However $\phi$ is also proximal, so we conclude
that $\phi$ is almost one-to-one. This shows that $\pi$ is an
HPI-extension contradicting our assumption.
\end{proof}

Thus combining this lemma with Proposition \ref{spi} we know
that $X_\infty$ contains a proximal chaotic subset $K$. However,
we do not know whether its image $\phi(K) \subset X$
is also such a set.

We conclude by formally stating our open problems.



\begin{prob}

1.\ 
A non-PI system contains a uniformly chaotic set
(Theorem \ref{nonpi}), is the converse true? (See remark \ref{w}.)

2.\ 
A strictly PI system which is not HPI contains a proximal chaotic subset
(Proposition \ref{spi}), is this true also for a PI non-HPI system?

\end{prob}
\medskip

\newpage

\section{Table}
In the table below we summarize the interrelations between
the various kinds of chaos discussed in the paper.
In each case the label refers to the existence of a large
chaotic set. We write `ch.' for chaos and `s' for strong. 
The labels ch. and s-ch. refer to the existence of Mycielski scrambled
set and Mycielski strongly scrambled set respectively.
Proximal means to say that the chaotic set in question is a proximal set.
\begin{table}[h]
\begin{center}

\begin{equation*}
\xymatrix
{
& {\rm weak\ mixing} \ar[d] & &\\
& {\rm dense\ uniform\ ch.} \ar[d]  &  & \\
{\rm Devaney\ ch.} \ar[r] & {\rm uniform\ ch.}\ar[d] & h_{\rm top}>0 \ar[d] &\\
& {\rm proximal\ s{\text -} ch.}\ar[r]\ar[d] & {\rm s{\text -}ch.} \ar[r]\ar[d] &
s{\text -}LY ch.\ar[d]\\
& {\rm proximal\ ch.}\ar[r] & {\rm ch.} \ar[r] & LY ch.
}
\end{equation*}

\medskip

\medskip

\caption{ \protect  Types of chaotic behavior}
\end{center}
\end{table}

\medskip


\section{Appendix}

\subsection{A characterization of RIC extensions}

Following usual notation we write $\beta\Z$ for the
\v{C}ech-Stone compactification of the integers,
and we fix a minimal left ideal $M \subset \beta\Z$
and an idempotent $u=u^2\in J(M)$, where $J(M)$ is the nonempty set of idempotents in $M$. Then the subset $G = uM$ is a maximal subgroup
of the semigroup $M$ which decomposes as a disjoint union
$M=\bigcup\{vG: v \in J(M)\}$. The group $G$ can be identified with
the group of automorphisms of the dynamical system $(M,\Z)$
(see e.g \cite{G76} or \cite{Au}). We also recall that the semigroup
$\beta\Z$ is a universal enveloping semigroup and thus ``acts" on
every compact $\Z$ dynamical system. In particular, when $(X,T)$ is
a dynamical system the homeomorphism $T$ defines in a natural way
a homeomorphism on $C(X)$, the compact space
of closed subsets of $X$.
Now for $p \in \beta\Z$ the ``action" of $p$ on
the point $A \in C(X)$ is well defined. In order to avoid confusion here we
denote the resulting element of $C(X)$ by $p \circ A$ and refer
to this action as the {\em circle operation}. A more concrete definition of
this set is
$$
p \circ A = \limsup T^{n_i}A,
$$
where, denoting by $S$ the generator of $\Z$, $S^{n_i}$ is any net in
$\beta \Z$ which converges to $p$. Thus we always have
$pA =\{px: x \in A\} \subset p \circ A$, but usually the inclusion is
proper, as often $pA$ is not even a closed subset of $X$.
A {\em quasifactor} of a system
$(X,T)$ is a closed invariant set $\mathcal{M} \subset C(X)$
such that $\bigcup \{A: A \in \mathcal{M}\} = X$.

Recall that an extension
$\pi : X \to Y$ of minimal dynamical systems is called a {\em relatively
incontractible (RIC) extension} if for every $p \in \beta \Z$
we have $p \circ Fx_0 = \pi^{-1}(py_0)$, where
$x_0 = ux_0$ is a point in $X$, $y_0=\pi(x_0)$ and 
$F = \mathfrak{G}(Y,y_0) = \{\alpha \in G: \alpha y_0 = y_0\}$ is the 
{\em Ellis group} of the pointed minimal system $(Y,y_0,T)$.

\begin{thm}\label{RIC}
The extension $\pi : X \to Y$ is RIC if and only if 
it is open and for every $n \ge 1$
the minimal points are dense in the relation
$$
R^n_\pi = \{(x_1,\dots,x_n) \in X^n : \pi(x_i)=\pi(x_j),\ \forall\   1\le i
\le j \le n\}.
$$
\end{thm}

\begin{proof}
Suppose $\pi$ is RIC. Then clearly  the map $y \mapsto \pi^{-1}(y)$
is continuous, i.e. $\pi$ is an open map.
Since every $n$-tuple $(x_1,\dots,x_n)$ with $x_i \in Fx_0,\ i=1,\dots,n$
is a minimal point of  $R^n_\pi$, the density of minimal
points in $R^n_\pi$ follows from the definition of
the circle operation.

Conversely, suppose $\pi$ is open and the minimal
points are dense in $R^n_\pi$.
We first note that the fact that $\pi$ is RIC
is equivalent to the statement that the minimal quasifactor
$\mathcal{M}=\{p \circ Fx_0 : p \in M\} \subset C(X)$,
where $M$ is a minimal ideal in $\beta \Z$,  coincides with the
collection $\{\pi^{-1}(y) : y \in Y\}$.
Thus it suffices to show that for an arbitrary $y \in Y$ the point
$\pi^{-1}(y)$ is in $\mathcal{M}$.

Let $d$ be any continuous semi-metric on $X$. Let $(x_1,\dots,x_n)$
be an $n$-tuple of elements of $\pi^{-1}(y)$ which is $\ep$-dense in
$\pi^{-1}(y)$ (with respect to  $d$). By assumption then, there is an $n$-tuple
$(x'_1,\dots,x'_n)$ of elements of $\pi^{-1}(y')$, such that
(i) $d(y,y')< \ep$, (ii) $d(x_i,x_i') < \ep$ for every $i$, and
(iii) $(x'_1,\dots,x'_n)$ is a minimal point of $R^n_\pi$.

There is a minimal idempotent $v \in J(M)$ such that
$vx'_i=x'_i$ for every $i$ and it follows that
$\{x'_1,\dots,x'_n\} \subset v \pi^{-1}(y')$.
Note that we must have $vy'=y'$ and there is therefore some
$p_{y'}= p = vp \in M$ with $y'=py_0$.

Now
$$
\{x'_1,\dots,x'_n\} \subset v \pi^{-1}(y') \subset v \circ v \pi^{-1}(y')
= v \circ pF x_0 = p \circ Fx_0 \subset \pi^{-1}(y').
$$
Since we have $\limsup_{y'\to y} \pi^{-1}(y') \subset \pi^{-1}(y)$, we
conclude that the set $\pi^{-1}(y)$ is a limit point
of the sets $p_{y'} \circ Fx_0 \in \mathcal{M}$ . Thus
also $\pi^{-1}(y)\in \mathcal{M}$ and the proof is complete.
\end{proof}

\medskip



\end{document}